\documentclass[twoside,6pt,leqno]{article}
\usepackage{textcase}
\usepackage{amsmath}
\usepackage{amssymb}
\usepackage{amsbsy}
\usepackage{amscd}
\usepackage{amsfonts}
\usepackage{amsopn}
\usepackage{amstext}
\usepackage{amsthm}
\usepackage{amsxtra}
\usepackage{graphicx,setspace}

\newtheorem{theorem}{Theorem}[section]
\newtheorem{lemma}{Lemma}[section]
\newtheorem{proposition}{Proposition}[section]
\newtheorem{corollary}{Corollary}[section]
\renewcommand{\theequation}{\mbox{\arabic{section}.\arabic{equation}}}
\def\reverseddots{\mathinner{\mkernlmu\raise7pt
\vbox{\kern7pt\hbox{.}}
\mkern-12mu\raise4pt\hbox{.}
\mkern-12mu\raise1pt\hbox{.}
\mkern1mu}}
\title{\bf{Bounds for $GL(3)\times GL(2)$ $L$-functions and $GL(3)$ $L$-functions}}
\author{Xiaoqing Li}
\date{}
\begin{document}
\maketitle
\pagenumbering{arabic}
\pagestyle{myheadings}
\markboth{}{}
\thispagestyle{headings}
\begin{abstract} In this paper, we will give the subconvexity bounds for self dual $GL(3)$ $L-$functions in the $t$ aspect as well as subconvexity bounds for 
self dual $GL(3)\times GL(2)$ $L-$functions in the $GL(2)$ spectral aspect.
\end{abstract}
\section
{Introduction}
Bounding $L$-functions on their critical lines is a far-reaching problem in number theory. For a general automorphic $L$-function, one may apply the Phragmen-Lindeloff interpolation method toghether with bounds on the $L$-function in $\Re s>1$ and $\Re s<0$ (the latter
coming from the functional equation) to give an upper bound for the $L$-function on the line $\Re s=\frac{1}{2}.$
The resulting bound is usually referred to as the convexity bound (or the trvial bound) for the $L$-function. While the Lindeloff hypothesis is still out of reach, breaking
the convexity bounds for $L$-functions is an interesting problem.\\
For $L$-functions of degree one, that is Dirichlet $L$-functions, such subconvexity estimates are due to Weyl [We] in the t-aspect and
 Burgess in the q-aspect[Bu]. For degree two $L$-functions this was achieved in a series of papers by Good [Go], Meurman
[Me] and especially Duke, Friedlander and Iwaniec [DFI1, DFI2, DFI3]. Subconvexity for Rankin-Selberg $L$-functions on $GL(2)\times GL(2)$ 
were known due to Sarnak [Sa], Kowalski, Michel and Vanderkam [KMV], Michel [Mi], Harcos and Michel [HM], Michel and Venkatesh [MV1], Lau, Liu and Ye 
[LYL], etc 
(see the references in [MV2]). Impressive subconvexity estimates for triple $L$-functions on $GL(2)$ were made by Bernstein and Reznikov [BR], see also Venkatesh [Ve].\\
Much less is known for subconvexity bounds for $L$-functions on higher rank groups. In this paper, we establish such subconvexity estimates for Rankin-Selberg $L$-functions on $GL(2)\times GL(3)$
and $L$-functions on $GL(3)$. To begin with, let $f(z)$ be a self dual Hecke-Maass form of type $(\nu, \nu)$ for $SL(3, \mathbb{Z}),$ normalized so that 
the first Fourier coefficient is $1.$ We define the $L$-function  
\begin{equation}
L(s, f)=\sum\limits_{m=1}^{\infty} A(m, 1)m^{-s}.
\end{equation}
For $f$ and each $u_j(z)$ in an orthonormal basis of even Hecke-Maass forms for $SL(2, \mathbb{Z}),$ we define the Rankin-Selberg $L$-function
\begin{equation}
L(s, f\times u_j)=
\sum\limits_{m\geqslant 1}\sum\limits_{n\geqslant 1}\frac{\lambda_j(n)A(n, m)}{(m^2n)^s}.
\end{equation} Our main theorem is the following:
\begin{theorem}
Let $f$ be a fixed self dual Hecke-Maass form for $SL(3, \mathbb{Z})$ and $u_j$ be an orthonormal basis
 of even Hecke-Maass forms for $SL(2, \mathbb{Z})$ corresponding to the Laplacian eigenvalue $\frac{1}{4}+t_j^2$ with 
$t_j\geqslant 0,$
 then for $\varepsilon>0, T$ large and
$T^{\frac{3}{8}+\varepsilon}\leqslant M\leqslant T^{\frac{1}{2}},$ 
we have
\begin{eqnarray}
\lefteqn{
{\sum_j}^{'} e^{-\frac{(t_j-T)^2}{M^2}}L\left(\frac{1}{2}, f\times u_j\right)
+
\frac{1}{4\pi}\int\limits_{-\infty}^\infty e^{-\frac{(t-T)^2}{M^2}}\left|L\left(\frac{1}{2}-it, f\right)\right|^2 dt}\\
&&\hspace{7.5cm}\ll_{\varepsilon, f}T^{1+\varepsilon}M\nonumber
\end{eqnarray}
where ' means summing over the orthonormal basis of even Hecke-Maass forms.
\end{theorem}
\noindent{\bf Remarks} 1. The second term in (1.3) comes from the Rankin-Selberg $L$-function of $f$ and the Eisenstein series on $GL(2).$\\
2. By considering the case that $f$ is the minimal Eisenstein series on $GL(3)$, one sees that the sign of the functional equation of 
$L(s, f\times u_j)$ is $+1$ when $u_j$ is an even Hecke-Maass form  and $-1$ when $u_j$ is an odd Hecke-Maass form for $SL(2, \mathbb{Z})$. For  this reason we restrict to even Hecke-Maass forms in (1.3). This feature doesn't appear if one averages the second moment of the 
$L$-functions.\\
3. Since $f$ is a self dual Hecke-Maass form  of $GL(3)$, it has to be orthogonal ([JS]) which means the (partial) 
$L$-function $L^S(s, f, sym^2)$ has a  pole at $s=1;$ since $u_j$ is a Maass form of $GL(2)$, it is symplectic which means $L^S(s, u_j, sym^2)$ has no pole at $s=1.$ Then Lapid's theorem [La] says that $L(\frac{1}{2}, f\times u_j)\geqslant 0.$ Due to this important property, we have
\begin{corollary}
Under the same assumptions as in the above theorem,
$$
L\left(\frac{1}{2}, f\times u_j\right)\ll_{\varepsilon, f}(1+|t_j|)^{\frac{11}{8}+\varepsilon}.
$$
\end{corollary}
The corresponding convexity bound for $L(\frac{1}{2}, f\times u_j)$ is 
$t_j^{\frac{3}{2}+\varepsilon}$ with $\varepsilon>0,$ so the above bound breaks the convexity
bound. \\
\noindent{\bf Remarks} 1. The nonnegativity of $L(\frac{1}{2}, f\times u_j)$ plays a crucial role in our approach. Otherwise, one can hardly motivate the goal of studying the first moment.\\
2. In the case that $f$ is an Eisenstein series on $GL(3)$, our approach recovers the subconvexity of a $GL(2)$ 
$L$-function in the eigenvalue aspect.\\
Ignoring the contribution of the cuspidal spectrum in (1.3) by the nonnegativity of $L(\frac{1}{2}, f\times u_j)$ [La], one has
$$
\int\limits_{-\infty}^\infty 
e^{-\frac{(t-T)^2}{T^{\frac{3}{4}}}}
\left|L\left(\frac{1}{2}-it, f\right)\right|^2 dt
\ll_{\varepsilon, f}T^{\frac{11}{8}+\varepsilon}.$$
By a standard argument [He], we have 
\begin{corollary}
For $f$ a self dual Hecke-Maass form for $SL(3, \mathbb{Z}),$
$$
L\left(\frac{1}{2}-it, f\right)
\ll_{\varepsilon, f}(|t|+1)^{\frac{11}{16}+\varepsilon}$$
where $\varepsilon>0.$
\end{corollary}
The corresponding convexity bound for $L(\frac{1}{2}-it, f)$ is $|t|^{\frac{3}{4}+\varepsilon}$ with $\varepsilon>0,$ so the above bound breaks the 
convexity bound for $L(\frac{1}{2}-it, f)$ in the $t$-aspect.
\noindent{\bf Remark.} Our method only breaks the convexity bounds of $L(\frac{1}{2}, f\times u_j)$  and $L(\frac{1}{2}, f)$ with $f$ self dual on $GL(3)$, i.e., $f$ comes from the symmetric lifts from $GL(2)$ (see [So]). New ideas are needed for the more general case $f$ is non self dual on $GL(3)$.\\
We end the introduction by a brief outline of the proof of the main theorem. Because we restrict to averaging over even Maass forms in (1.3), applying the approximate functional equation for the Rankin-Selberg $L$-functions and Kuznetsov's formula leads to two parts: 
$ {\tilde{\cal R}}_3^+$ (see (4.17)) -- weighted sums of Kloosterman sums twisted by 
$e^{\frac{4\pi i\sqrt{n}}{c}}$ and $\tilde{\cal R}_{2}^-$ (see (5.10)) -- weighted sums of Kloosterman sums without twisting. Instead of using Weil's bound for the Kloosterman sum which only leads to the convexity bound for the individual 
$L$-function, we expand the Kloosterman sums and makes crucial use of the Voronoi formula on $GL(3)$.  
$\tilde{\cal R}_{2}^-$ involves no twisting which allows a direct application of the Voronoi formula. 
$\tilde{\cal R}_3^+$ seems harder. However, as a miracle, the application of the Voronoi formula to 
$\tilde{\cal R}_3^+$ brings the twists by $e^{\frac{4\pi i\sqrt{n}}{c}}$ to twists by additive characters (see (4.24)). This breaks the duality of the Voronoi formula. A second application of the Voronoi formula twisted by additive characters then completes the estimation of
 ${\cal R}_3^+$. In using the Voronoi formula, one needs the asymptotic behavior of the integral transformations of the test functions. This is provided in Lemma 2.1. In the appendix, suggested by Sarnak, we also considered the subconvexity of the Rankin-Selberg
 $L$-function $L(s, f\times h)$ where $f$ is self dual on $GL(3)$ and $h$ runs through holomorphic forms of weight $k$ congruent to $0$ modulo $4$. 
The analysis is essentially the same as the nonholomorphic case.\\

The Voronoi formula for $GL(3)$ was first derived by Miller and Schmidt [MS] (see [GL] for a simple proof). It was first used by Sarnak and Watson to prove a Lindeloff like bound for the $L^4$ norm of a Maass form for $GL(2).$ For other
applications, see [Mi] and [Li]. Throughout the paper, $e(x)$ means $e^{2\pi ix}$ and negligible means $O(T^{-A})$ for any $A>0.$
\section {A review of automorphic forms}
\setcounter{equation}{0} 
In this section, we introduce notations and recall some standard facts of Maass forms for $GL(2)$ and $GL(3).$ We start from the upper half plane $\mathbb{H}.$
The Laplace operator
$$
\Delta=-y^2\left(
\frac{\partial^2}{\partial x^2}+\frac{\partial^2}{\partial y^2}\right)
$$
 has a spectral decomposition on $L^2(SL(2, \mathbb{Z})\setminus\mathbb{H}):$
$$
L^2(SL(2, \mathbb{Z})\setminus\mathbb{H})={\cal{C}}\oplus C(SL(2, \mathbb{Z})\setminus{\mathbb{H}})\oplus
{\cal{E}}(SL(2, \mathbb{Z})\setminus\mathbb{H}).
$$
Here $\cal {C}$ is the space of constant functions. $C(SL(2, \mathbb{Z})\setminus{\mathbb{H}})$ is the space of Maass forms 
and ${\cal{E}}(SL(2, \mathbb{Z})\setminus\mathbb{H})
$ is the space of Eisenstein series.\\
Let ${\cal{U}}=\{{u_j: j\geqslant 1}\}$ be an orthonormal basis of Hecke-Maass forms corresponding to the Laplacian eigenvalue $\frac{1}{4}+t_j^2$ 
with $t_j\geqslant 0$ in the 
space
$C(SL(2, \mathbb{Z})\setminus\mathbb{H}).$ Any $u_j(z)$ has the Fourier expansion
$$
u_j(z)=\sum\limits_{n\neq 0}\rho_j(n)W_{s_j}(nz)
$$
where $W_s(z)$ is the Whittaker function given by
$$
W_s(z)=2|y|^{\frac{1}{2}}K_{s-\frac{1}{2}}(2\pi|y|)e(x)
$$
and $K_s(y)$ is the $K$-Bessel function with $s=\frac{1}{2}+it$. $C(SL(2, \mathbb{Z})\setminus\mathbb{H})$ consists of even Maass forms and odd Maass forms
according to $u_j(-\bar{z})=u_j(z)$ or $u_j(-\bar{z})=-u_j(z).$ We can assume $u_j$ are eigenfunctions of all the Hecke operators corresponding to
the Hecke eigenvalue $\lambda_j(n).$ Then we have the formula
$$
\rho_j(\pm n)=\rho_j(\pm 1)\lambda_j(n)n^{-\frac{1}{2}}
$$
if $n>0.$ The Eisenstein series $E(z, s)$ defined by
\begin{equation}
E(z, s)=\frac{1}{2}\sum_{\substack{c, d\in \mathbb{Z}\\ (c, d)=1}}\frac{y^s}{|cz+d|^{2s}}
\end{equation}
has the following Fourier expansion
$$
E(z, s)=y^s+\phi(s)y^{1-s}+\sum\limits_{n\neq 0}\phi(n, s)W_s(nz)
$$
where
$$
\phi(s)=\sqrt{\pi}\frac{\Gamma(s-\frac{1}{2})}{\Gamma(s)}\frac{\zeta(2s-1)}{\zeta(2s)}
$$ with $\zeta(s)$ be the Riemann zeta function
and
$$
\phi(n, s)=\pi^s\Gamma(s)^{-1}\zeta(2s)^{-1}|n|^{-\frac{1}{2}}\eta(n, s)
$$
with
\begin{equation}
\eta(n, s)=\sum\limits_{ad=|n|}\left(\frac{a}{d}\right)^{s-\frac{1}{2}}.
\end{equation}
For any $m, n\geqslant 1$ and any test function $h(t)$ which is even and satisfies the following conditions:
\\
I) $h(t)$ is holomorphic in $|\Im t|\leqslant \frac{1}{2}+\varepsilon;$\\
II)$h(t)\ll (|t|+1)^{-2-\varepsilon}$ in the above strip,
we have the following Kuznetsov formula
(see [CI])
\begin{eqnarray}
\lefteqn{\;\;\;\;\;\;\;\;\;
{\sum\limits_{j\geqslant 1}}^{'}h(t_j)\omega_j\lambda_j(m)\lambda_j(n)\!+\!\frac{1}{4\pi}\!\!
\int\limits_{-\infty}^{\infty}\!
h(t)\omega(t)
\bar{\eta}\left(m, \frac{1}{2}+it\right)\eta\left(n, \frac{1}{2}+it\right)\!dt}\\
&&\!\!\!\!\!\!\!\!\!\!=\frac{1}{2}\delta(m, n)H
\!+\!\!\sum\limits_{c>0}\frac{1}{2c}\left\{S(m, n; c)H^+\!\!\left(\frac{4\pi\sqrt{mn}}{c}\right)\!+\!
 S(-m, n; c)H^-\!\!\left(\frac{4\pi\sqrt{mn}}{c}\right)\right\}\nonumber
\end{eqnarray}
where ${\sum}^{'}$ restricts to the even Maass forms, $\delta(m, n)$ is the Kronecker symbol,
$$
\omega_j=4\pi|\rho_j(1)|^2/\cosh \pi t_j,
$$
$$
\omega(t)=4\pi\Bigg|\phi\left(1, \frac{1}{2}+it\right)\Bigg|^2\cosh^{-1}\pi t,
$$
$$
H=\frac{2}{\pi}\int\limits_{0}^{\infty}h(t)\tanh (\pi t) tdt,
$$
$$
H^+(x)=2i\int\limits_{-\infty}^{\infty}J_{2it}(x)\frac{h(t)t}{\cosh\pi t}dt,
$$
$$
H^-(x)=\frac{4}{\pi}\int\limits_{-\infty}^\infty K_{2it}(x)\sinh (\pi t)h(t)tdt,
$$
$$
S(a, b; c)=\displaystyle{\sum_{d\bar {d}\equiv1(\text{mod}\; c)}}e\Big(\frac{da+\bar{d}b}{c}\Big)
$$
is the classical Kloosterman sum, in the above, $J_\nu(x)$ and $K_{\nu}(x)$ are the standard  $J-$Bessel function and $K-$Bessel function respectively.\\

Now we recall some background on Maass forms for $GL(3).$ We will follow the notations in Goldfeld's book [Gol]. Let $f$ be a Maass form of type $\nu=(\nu_1, \nu_2)$ for 
$SL(3, \mathbb{Z}).$ Thanks to Jacquet, Piatetskii-Shapiro and Shalika, we have the following Fourier Whittaker expansion 
\begin{equation}
f(z)=
\!\!\!\sum\limits_{\gamma\in U_{2}(\mathbb{Z})\backslash SL(2, \mathbb{Z})}
\;\sum\limits_{m_1=1}^{\infty}
\sum\limits_{m_{2}\neq 0}
\frac{A(m_1,  m_{2})}{m_1|m_2|}
W_{\rm J}\left(M\left(\!\!\begin{array}{rr}
                   \gamma\!&\! \\
                   \!&\! 1\\
                   \end{array}\!\!\right)z, \nu, \psi_{1,1}\right)
\end{equation}
where $U_2(\mathbb Z)$ is the group of $2\times 2$ upper triangular matrices with integer entries and ones on the diagonal, $ W_{\rm J}(z, \nu, \psi_{1, 1})$ is the Jacquet-Whittaker function and 
$M=\mbox{diag}\left(m_1|m_2|,m_1, 1\right).$ Set 
$$
\alpha=-\nu_1-2\nu_2+1,\;\; \beta=-\nu_1+\nu_2,\;\; \gamma=2\nu_1+\nu_2-1,
$$
for $k=0, 1;$ for $\psi(x)$ a smooth compactly supported function on $(0, \infty)$ and 
$\tilde{\psi}(s):=\int\limits_0^\infty\psi(x)x^s\frac{dx}{x},$
set
\begin{equation}
\Psi_k(x):=\int\limits_{\Re s=\sigma}(\pi^3x)^{-s}\frac{
\Gamma\left(\frac{1+s+2k+\alpha}{2}\right)
\Gamma\left(\frac{1+s+2k+\beta}{2}\right)
\Gamma\left(\frac{1+s+2k+\gamma}{2}\right)}
{\Gamma\left(\frac{-s-\alpha}{2}\right)
\Gamma\left(\frac{-s-\beta}{2}\right)
\Gamma\left(\frac{-s-\gamma}{2}\right)}
\tilde{\psi}(-s-k)ds
\end{equation}
with $\sigma>\max\{-1-\Re\alpha, -1-\Re\beta, -1-\Re\gamma\},$
$$
\Psi_{0, 1}^0(x)=\Psi_0(x)+\frac{\pi^{-3}c^3m}{n_1^2n_2i}\Psi_1(x)
$$
and
$$
\Psi_{0, 1}^1(x)=\Psi_0(x)-\frac{\pi^{-3}c^3m}{n_1^2n_2i}\Psi_1(x),
$$
 we have the following Voronoi formula on ${\rm GL}(3):$

\begin{proposition} {\rm{([MS], [GL])}} Let $\psi(x)\in C_c^{\infty}(0, \infty).\;$Let $A(m, n)$ denote the 
$(m, n)$-th Fourier coefficient of a Maass form for $SL(3, \mathbb{Z})$ as in 
(2.4).  Let
${d, \bar{d}, c} \in\mathbb{Z}$ with $c\neq 0, (d, c)=1,$
and $d\bar{d}\equiv 1 (\rm {mod} \;c).$ Then
we have
\begin{eqnarray}
\lefteqn{
\sum\limits_{n>0}A(m, n)e\left(\frac{n\bar{d}}{c}\right)\psi(n)\nonumber}\\
&&
=\frac{c\pi^{-\frac{5}{2}}}{4i}\sum\limits_{n_1|cm}\sum\limits
_{n_2>0}\frac{A(n_2, n_1)}{n_1n_2}S(md, n_2; m cn_1^{-1})\Psi_{0, 1}^0
\left(\frac{n_2n_1^2}{c^3m}\right)
\nonumber\\
&&\;\;\;\;+\frac{c\pi^{-\frac{5}{2}}}{4i}\sum\limits_{n_1|cm}\sum\limits
_{n_2>0}\frac{A(n_2, n_1)}{n_1n_2}S(m d, -n_2; m cn_1^{-1})
\Psi_{0, 1}^1
\left(\frac{n_2n_1^2}{c^3m}\right)
\nonumber,
\end{eqnarray}
where $S(a, b; c)$ is the Kloosterman sum defined as the above.
\end{proposition}
To apply Proposition 2.1 in practice, one needs to know the asymptotic behaviour of $\Psi_0(x)$ and $\Psi_1(x).$ By changing variables $s+1\rightarrow s$ in the definition of $\Psi_1(x)$, one sees that $x^{-1}\Psi_1(x)$ has similar asymptotic behavior as of $\Psi_0(x)$. Therefore, in the following, we only consider $\Psi_0(x)$.
\begin{lemma}{\rm{([Li])}}
Suppose $\psi(x)$ is a smooth function compactly supported on $[X, 2X],$ $\Psi_0(x)$ is defined by (2.5), then for any fixed integer $K\geqslant 1$ and
$xX\gg 1,$ we have
\begin{eqnarray}
\lefteqn{
\Psi_0(x)=2\pi^4xi\int\limits_0^\infty\psi(y)\sum\limits_{j=1}^K
\frac{c_j\cos (6\pi x^{\frac{1}{3}}y^{\frac{1}{3}})+d_j\sin (6\pi 
x^{\frac{1}{3}}y^{\frac{1}{3}})}{(\pi^3xy)^{\frac{j}{3}}}dy}\nonumber\\
&&\hspace{6cm}+O\left((xX)^{\frac{-K+2}{3}}\right),
\nonumber
\end{eqnarray}
where $c_j$ and $d_j$ are constants depending on $\alpha, \beta$ and $\gamma,$ in particular, \\$c_1=0, d_1=-\frac{2}{\sqrt{3\pi}}.$
\end{lemma} 
\noindent{\bf Remark.} When $xX\ll 1,$ moving the line of integration to $\sigma=-\frac{11}{20}$, by Stirling's formula for the $\Gamma$ functions and integration by
part once for $\tilde{\psi}(s),$ one shows that
$$
\Psi_0(x)\ll \int\limits_0^\infty |\psi^{'}(x)|dx.
$$
Note that a special case of the above lemma (when $\alpha=\beta=\gamma=0$ ) was given by Ivic (see [Iv]).
Now let $f$ be a self dual Hecke-Maass form of type $(\nu, \nu)$ for $SL(3, \mathbb{Z}),$ normalized to have the first Fourier coefficient $A(1, 1)$ equal to $1.$
We associate the $L$-function $L(s, f)$ defined by (1.1). It is entire and satisfies
the functional equation
$$
G_{\nu}(s)L(s, f)=G_{\nu}(1-s)L(1-s, f)
$$
where
$$
G_{\nu}(s)=\pi^{\frac{-3s}{2}}
\Gamma\left(\frac{s+1-3\nu}{2}\right)
\Gamma\left(\frac{s}{2}\right)
\Gamma\left(\frac{s-1+3\nu}{2}\right).
$$
The Rankin-Selberg $L$-function defined by
$$
L(s, f\times f):=\sum\limits_{m\geqslant 1}\sum\limits_{n\geqslant 1}\frac{|A(m, n)|^2}{(m^2n)^s}
$$
for $\Re s$ large has a meromorphic continuation to the whole plane with the only simple pole at $s=1.$ By a standard contour integration,
one shows that
\begin{equation}
\mathop{\sum\sum}_{m^2n\leqslant N}
|A(m, n)|^2\ll_f N.
\end{equation}
By Cauchy's inequality and (2.6), one derives that
\begin{equation}
\sum\limits_{n\leqslant N}|A(m, n)|\ll_f N|m|.
\end{equation}
The Rankin-Selberg $L$-function of $f$ and $u_j$
defined by (1.2) is entire and satisfies the functional equation
\begin{equation}
\Lambda(s, f\times u_j)=\Lambda(1-s, f\times u_j)
\end{equation}
where
\begin{eqnarray}
\lefteqn{
\Lambda(s, f\times u_j)=\pi^{-3s}
\Gamma\left(\frac{s-it_j-\alpha}{2}\right)
\Gamma\left(\frac{s-it_j-\beta}{2}\right)
\Gamma\left(\frac{s-it_j-\gamma}{2}\right)}\nonumber\\
&&\hspace{1.5cm}\times
\Gamma\left(\frac{s+it_j-\alpha}{2}\right)
\Gamma\left(\frac{s+it_j-\beta}{2}\right)
\Gamma\left(\frac{s+it_j-\gamma}{2}\right)L(s, f\times u_j)
\nonumber
\end{eqnarray}
and 
\begin{equation}
\alpha=-3\nu+1,\;\;\beta=0,\;\;\gamma=3\nu-1.
\end{equation}
To the above Maass form $f$ and the Eisenstein series $E\left(z, \frac{1}{2}+it\right)$ (recall (2.1))
we associate the $L$-function
$$
L(s, f\times E):=\sum\limits_{m\geqslant 1}\sum\limits_{n\geqslant 1}\frac{\bar{\eta}(n, \frac{1}{2}+it)A(n, m)}{(m^2n)^s}
.$$
By looking at the Euler products
$$
L(s, f)=\sum\limits_{n\geqslant 1}\frac{A(n, 1)}{n^s}=
\prod\limits_p\prod\limits_{i=1}^3(1-\beta_{p, i}p^{-s})^{-1},
$$
$$
L(s, E)=\sum\limits_{n\geqslant 1}\eta\left(n, \frac{1}{2}+it\right)n^{-s}=\prod\limits_p(1-p^{-s+it})^{-1}(1-p^{-s-it})^{-1},
$$
one derives that (see [Gol] pp. 379)
\begin{eqnarray}
\lefteqn{
L(s, f\times E)=\prod\limits_p\prod\limits_{k=1}^3
(1-\beta_{p, k}p^{it-s})^{-1}(1-\beta_{p, k}p^{-it-s})^{-1}}\nonumber\\
&&\hspace{5cm}=L(s-it, f)L(s+it, f).\nonumber
\end{eqnarray}
It yields that 
$$
L\left(\frac{1}{2}, f\times E\right)=\left|L\left(\frac{1}{2}-it, f\right)\right|^2.
$$
This satisfies the functional equation (2.8) which can also be verified directly using the functional equation of $L(s, f).$
Set
$$
F(u)=\left(\cos \frac{\pi u}{A}\right)^{-3A},
$$ for $|\Im t|\leqslant 1000,$ where $A$ is a positive integer, 
\begin{equation}
V(y, t)=\frac{1}{2\pi i}\int\limits_{(1000)}y^{-u}F(u)\frac{\gamma(\frac{1}{2}+u, t)}{\gamma(\frac{1}{2}, t)}
\frac{du}{u}
\end{equation} and 
\begin{eqnarray}
\lefteqn{\;\;\;\;\;\;\;\;\;\;\;
\gamma(s, t)=\pi^{-3s}
\Gamma\left(\frac{s-it-\alpha}{2}\right)
\Gamma\left(\frac{s-it-\beta}{2}\right)
\Gamma\left(\frac{s-it-\gamma}{2}\right)}\nonumber\\
&&\hspace{3.9cm}\times
\Gamma\left(\frac{s+it-\alpha}{2}\right)
\Gamma\left(\frac{s+it-\beta}{2}\right)
\Gamma\left(\frac{s+it-\gamma}{2}\right).
\nonumber
\end{eqnarray}
The integral is justified by Luo-Rudnick-Sarnak's 
bound on the Ramanujan conjecture $|\Re \alpha|, |\Re\beta|, |\Re \gamma|\leqslant \frac{1}{2}-\frac{1}{10}$ (see [LRS]). One has the following 
approximate functional equation for $L(s, f\times u_j)$ (see [IK] or [Li]):
\begin{lemma} 
For a self dual Maass form $f$ of type $(\nu, \nu)$ for $SL(3, \mathbb{Z})$ and any $u_j(z)$ associated to the Laplacian eigenvalue $\frac{1}{4}+t_j^2$ in the orthonormal
basis of even Hecke-Maass forms for $SL(2, \mathbb{Z}),$ 
we have
\begin{equation}
L\left(\frac{1}{2}, f\times u_j\right)=2\sum\limits_{m\geqslant 1}\sum\limits_{n\geqslant 1}
\frac{\lambda_j(n)A(n, m)}{(m^2n)^{\frac{1}{2}}}
V(m^2n, t_j).
\end{equation} 
\end{lemma}
$V(y, t)$ has the following properties which effectively limit the terms in (2.11) with $m^2n\ll |t_j|^3.$
\begin{lemma}For $y, t>0, i=1, 2,$\\
{\rm 1)} the derivatives of $V(y, t)$ with respect to $y$ satisfy 
$$
y^a\frac{\partial^a}{\partial y^a}V(y, t)\ll \left(1+\frac{y}{|t|^3}\right)^{-A},
$$
$$
y^a\frac{\partial^a}{\partial y^a}V(y, t)=\delta_a+O\left(\left(\frac{y}{|t|^3}\right)^c\right),
$$
where $0<c\leqslant \frac{1}{3}\text{min} \{
\frac{1}{2}-\Re \alpha, \frac{1}{2}-\Re \beta, \frac{1}{2}-\Re \gamma\},\;
\delta_0=1, 0$ otherwise and the implied constants depend only on $c, a, A$, $\alpha, \beta$ and $\gamma$.\\
{\rm 2)} if $1\leqslant y\ll t^{3+\varepsilon},$ then as $t\rightarrow\infty,$ we have
\begin{eqnarray}
\lefteqn{
V(y, t)=\frac{1}{2\pi i}
\int\limits_{\left(\frac{1}{2}\right)}\!\!
\left(\frac{t^3}{8\pi^3 y}\right)^uF(u)\left[1+\frac{p_1(v)}{t}
+\cdots+\frac{p_{n-1}(v)}{t^{n-1}}+O\left(\frac{p_n(v)}{t^n}\right)\right]\frac{du}{u}\nonumber}\\
&\hspace{10.2cm}+O\left(t^{-B}\right)\nonumber
\end{eqnarray}
where $v=\Im u,$ $p_i(v)$ are polynomials of $v$ and $B$ is arbitrarily large.
\end{lemma}
\noindent{\bf Proof.} 1) See [IK], pp. 100.\\
2) It follows from Stirling's formula
$$
\log \Gamma(s+b)=\left(s+b-\frac{1}{2}\right)\log s-s+\frac{1}{2}\log 2\pi+\sum\limits_{j=1}^k\frac{a_j}{s^j}+O_\delta\left(\frac{1}{|s|^{k+1}}\right),
$$
which is valid for $b$ a constant, any fixed integer $K\geqslant 1, |\arg s|\leqslant \pi-\delta$ for $\delta>0,$ where the point $s=0$ and the neighbourhoods of the poles of $\Gamma(s+b)$ are excluded,
and the $a_j$ are suitable constants.  \mbox{$\Box$}\\
$L(s, f\times E)$ has the similar approximate functional equation as the above
\begin{equation}
L\left(\frac{1}{2}, f\times E\right)=2\sum\limits_{m\geqslant 1}\sum\limits_{n\geqslant 1}
\frac{\eta(n, \frac{1}{2}+it)A(n, m)}{(m^2n)^{\frac{1}{2}}}
V(m^2n, t).
\end{equation} 
Now we introduce the spectrally normalized first moment of the central values of $L$-functions
\begin{equation}
W:={\sum_j}^{'} e^{-\frac{(t_j-T)^2}{M^2}}\omega_jL\left(\frac{1}{2}, f\times u_j\right)
+
\frac{1}{4\pi}\int\limits_{-\infty}^\infty e^{-\frac{(t-T)^2}{M^2}}\omega(t)\left|L\left(\frac{1}{2}-it, f\right)
\right|^2 dt
\end{equation}
where $\omega_j$ and $\omega(t)$ are defined below (2.3). Due to Iwaniec [Iw2], we know
$$
\omega_j\gg t_j^{-\varepsilon}
$$
and as a well-known fact ([Ti], pp. 111) we also know 
$$
\omega(t)\gg t^{-\varepsilon},
$$
one has
$$
{\sum_j}^{'} e^{-\frac{(t_j-T)^2}{M^2}}L\left(\frac{1}{2}, f\times u_j\right)
+
\frac{1}{4\pi}\int\limits_{-\infty}^\infty e^{\frac{-(t-T)^2}{M^2}}\left|L\left(\frac{1}{2}-it, f\right)\right|^2 dt\ll WT^{\varepsilon}
$$
for any $\varepsilon>0.$ Therefore, for Theorem 1.1 we need to show that
\begin{equation}
W\ll_{\varepsilon, f}T^{1+\varepsilon} M.
\end{equation}
To use the Kuznetsov formula, the test function has to be even. For that purpose, we introduce
\begin{equation}
{\cal{W}}:={\sum_j}^{'} k(t_j)\omega_jL\left(\frac{1}{2}, f\times u_j\right)
+
\frac{1}{4\pi}\int\limits_{-\infty}^\infty k(t)\omega(t)\left|L\left(\frac{1}{2}-it, f\right)\right|^2 dt,
\end{equation}
here
\begin{equation}
k(t)=e^{-\frac{(t-T)^2}{M^2}}+ e^{-\frac{(t+T)^2}{M^2}}.
\end{equation}
Applying (2.11) and (2.12) to $\cal{W},$ by smooth dyadic subdivisions it suffices for our purposes to estimate sums of the form
\begin{eqnarray}
\lefteqn{\;\;\;\;\;\;\;\;\;\;\;{\cal R:=}
2{\sum\limits_j}^{'}k(t_j)\omega_j\sum\limits_{m\geqslant 1}\sum\limits_{n\geqslant 1}
\frac{\lambda_j(n)A(n, m)}
{(m^2n)^{\frac{1}{2}}}V(m^2n, t_j)g\left(\frac{m^2n}{N}\right)}\\
&&\;\;\;\;\;\;\;\;\;\;+\frac{2}{4\pi}\int\limits_{-\infty}^\infty k(t)\omega(t)\sum\limits_{m\geqslant 1}\sum\limits_{n\geqslant 1}
\frac{\eta(n, \frac{1}{2}+it)A(n, m)}
{(m^2n)^{\frac{1}{2}}}V(m^2n, t)g\left(\frac{m^2n}{N}\right)dt\nonumber.
\end{eqnarray}
Here $g$ is essentially a fixed smooth function of compact support on $[1, 2]$ and $N$ is at most $T^{3+\varepsilon}, \varepsilon>0.$
We then transform $\cal {R}$ by the Kuznetsov formula (2.3) into
\begin{equation}
{\cal R}={\cal D}+{\cal R}^+ +{\cal R}^-
\end{equation}
where
\begin{equation}
{\cal D}=
\sum\limits_{m\geqslant 1}\sum\limits_{n\geqslant 1}
\frac{A(n, m)}{(m^2n)^{\frac{1}{2}}}g\left(\frac{m^2n}{N}\right)\delta(n, 1)H_{m, n}
\end{equation}
is the contribution of the diagonal term with 
\begin{equation}
H_{m, n}=\frac{2}{\pi}\int\limits_0^\infty k(t)V(m^2n, t)\tanh( \pi t )tdt,
\end{equation}
\begin{equation}
{\cal R}^+ = \sum\limits_{m\geqslant 1}\sum\limits_{n\geqslant 1}\frac{A(n, m)}{(m^2n)^{\frac{1}{2}}}g\left(\frac{m^2n}{N}\right)
\sum\limits_{c>0}c^{-1}S(n, 1; c)H_{m, n}^+\left(\frac{4\pi\sqrt{n}}{c}\right)
\end{equation}
with
\begin{equation}
H_{m, n}^+(x)=2i\int_{-\infty}^{\infty} J_{2it} (x)\frac{k(t)V(m^2n, t)t}{\cosh \pi t}dt
\end{equation}
and 
\begin{equation}
{\cal R}^- = \sum\limits_{m\geqslant 1}\sum\limits_{n\geqslant 1}\frac{A(n, m)}{(m^2n)^{\frac{1}{2}}}g\left(\frac{m^2n}{N}\right)
\sum\limits_{c>0}c^{-1}S(n, 1; c)H_{m, n}^-\left(\frac{4\pi\sqrt{n}}{c}\right)
\end{equation} with
\begin{equation}
H_{m, n}^-(x)=\frac{4}{\pi}\int_{-\infty}^\infty K_{2it}(x)\sinh (\pi t) k(t)V(m^2n, t)tdt.
\end{equation}
The next three sections are devoted to the estimation of ${\cal D}, {\cal R}^+,$ and  ${\cal R}^-$ respectively. 
\section{The diagonal terms}
\setcounter{equation}{0} 
Recall that $\cal {D}$ is the contribution to $\cal {R}$ (see (2.18)) from the diagonal terms defined by (2.19). Obviously
$$
{\cal D}=\sum\limits_{m\geqslant 1}\frac{A(1, m)}{m}g\left(\frac{m^2}{N}\right)H_{m, 1}
$$
where
\begin{eqnarray}
\lefteqn{
H_{m, 1}=
\frac{2}{\pi}\int\limits_0^\infty \left[e^{-\frac{(t-T)^2}{M^2}}+e^{-\frac{(t+T)^2}{M^2}}\right]
V(m^2, t)\tanh (\pi t) tdt}\\
&&\hspace{0.25cm}=\frac{2}{\pi}\int\limits_0^\infty e^{-\frac{(t-T)^2}{M^2}}
V(m^2, t)\tanh (\pi t) tdt + O(T^{-A})\nonumber
\end{eqnarray}
with $A$ arbitrarily large.
By Lemma 2.3 and (2.7), we have
\begin{equation}
\sum\limits_{m\geqslant 1}\frac{A(1, m)}{m}g\left(\frac{m^2}{N}\right)V(m^2, t)\ll_{\varepsilon, f}
\left(|t|+1\right)^{\varepsilon}.
\end{equation}
It follows from (3.1) and (3.2) that
$$
{\cal D}\ll_{\varepsilon, f}T^{1+\varepsilon}M
$$ as we want.
\section{The terms related to the $J-$Bessel function}
\setcounter{equation}{0} 
This section is devoted to the estimation of ${\cal R}^+$ which is defined by (2.21).\\
We split ${\cal R}^+$ into three parts ${\cal R}^+_1, {\cal R}^+_2, {\cal R}^+_3 $ with
\begin{equation}
{\cal R}^+_1 = \sum\limits_{m\geqslant 1}\sum\limits_{n\geqslant 1}\frac{A(n, m)}{(m^2n)^{\frac{1}{2}}}g\left(\frac{m^2n}{N}\right)
\sum\limits_{c\geqslant C_1/m}c^{-1}S(n, 1; c)H_{m, n}^+\left(\frac{4\pi\sqrt{n}}{c}\right),
\end{equation}
\begin{equation}
{\cal R}^{+}_2 = \sum\limits_{m\geqslant 1}\sum\limits_{n\geqslant 1}\frac{A(n, m)}{(m^2n)^{\frac{1}{2}}}g\left(\frac{m^2n}{N}\right)
\sum\limits_{C_2/m\leqslant c\leqslant C_1/m }c^{-1}S(n, 1; c)H_{m, n}^+\left(\frac{4\pi\sqrt{n}}{c}\right),
\end{equation}
\begin{equation}
{\cal R}^+_3 = \sum\limits_{m\geqslant 1}\sum\limits_{n\geqslant 1}\frac{A(n, m)}{(m^2n)^{\frac{1}{2}}}g\left(\frac{m^2n}{N}\right)
\sum\limits_{c\leqslant C_2/m}c^{-1}S(n, 1; c)H_{m, n}^+\left(\frac{4\pi\sqrt{n}}{c}\right)
\end{equation}
where
\begin{equation}
C_1=T, \;\;\;\;C_2=\frac{\sqrt{N}}{T^{1-\varepsilon}M}.
\end{equation}
First we will estimate (4.1). Recall $H_{m, n}^+(x)$ is defined by (2.22). Moving the line of integration to $\Im t=-100,$ $H_{m, n}^+(x)$ becomes
 \begin{equation}
2i\!\!\int\limits_{-\infty}^{\infty}\!J_{2iy+200}( x)\frac{k(-100i+y)V(m^2n, -100 i+y)(-100 i+y)}
{\cosh \pi (-100 i+y)}
dy.
\end{equation}
By the integral representation of the $J-$Bessel function ([GR], 8.411 4)
$$
J_\nu(z)=2\frac{\left(\frac{z}{2}\right)^\nu}{\Gamma(\nu+\frac{1}{2})\Gamma(\frac{1}{2})}\int\limits_0^{\frac{\pi}{2}}\sin^{2\nu}\theta\cos(z\cos\theta)d\theta
$$
for $\Re \nu>-\frac{1}{2},$ one derives that
\begin{equation}
J_{2iy+200}(x)\ll\left(\frac{x}{|y|}\right)^{200} e^{\pi|y|}.
\end{equation}
Using Stirling's formula, we have
\begin{equation}
V(m^2n, -100i+y)\ll\left(\frac{|y|^3}{m^2n}\right)^{100}.
\end{equation}
Combining (4.5), (4.6) and (4.7), we have
\begin{equation}
H_{m, n}^+(x)\ll x^{200}T^{100}(m^2n)^{-100}TM.
\end{equation}
Thus, by (2.7), (4.8) and the trivial bound for the Kloosterman sum, one concludes that
\begin{equation}
{\cal R}^+_1\ll N^{\frac{1}{2}}T^{-98}M\ll 1.
\end{equation}
Next we will estimate ${\cal R}^+_2.$ By [GR] (8.411 11), one derives that
$$
\frac{J_{2it}( x)+J_{-2it}(x)}{\cosh \pi t}=\frac{2}{\pi}\int\limits_{-\infty}^{\infty}
\sin ( x\cosh \zeta)e\left( \frac{t\zeta}{\pi}\right) d\zeta.
$$
Applying the above integral representation and partial integration in $\zeta$ once, we have
\begin{eqnarray}
\lefteqn{
H_{m, n}^+(x)=
\frac{4i}{\pi}\int\limits_{t=0}^\infty\int\limits_{\zeta=-T^\varepsilon}^{T^\varepsilon}
 te^{-\frac{(t-T)^2}{M^2}}V(m^2n, t)\sin(x\cosh\zeta)e\left(\frac{t\zeta}{\pi}\right)dtd\zeta}\nonumber\\
&&\hspace{9cm}+O(T^{-A})\nonumber
\end{eqnarray}
with $A$ arbitrarily large. By changing variables $\frac{t-T}{M}\rightarrow t,$ we have
\begin{eqnarray}
\lefteqn{
H_
{m, n}^+(x)=
\frac{4iM}{\pi}\int\limits_{t=-\frac{T}{M}}^\infty\int\limits_
{\zeta=-T^\varepsilon}^{T^\varepsilon}
 (T+tM)e^{-t^2}V(m^2n, tM+T)\sin(x\cosh\zeta)  }\nonumber\\
&&\hspace{6cm}\times e\left(\frac{(tM+T)\zeta}{\pi}\right)dtd\zeta+O(T^{-A}).\nonumber
\end{eqnarray}
Extending the $t$ integral to $(-\infty, \infty)$ with a negligible error term, we have
$$
H_{m, n}^+(x)=H_{m, n}^{+, 1}(x) +H_{m, n}^{+, 2}(x) + O(T^{-A})
$$
where
\begin{eqnarray}
\lefteqn{
H_
{m, n}^{+, 1}(x)=
\frac{4iMT}{\pi}\int\limits_
{t=-\infty}^\infty
\int\limits_{\zeta=- T^\varepsilon}^{\zeta=T^\varepsilon}
 e^{-t^2}V(m^2n, tM+T)\sin(x\cosh\zeta)e\left(\frac{tM\zeta}{\pi}\right)\nonumber}\\
&&\hspace{8.9cm}\times e\left(\frac{T\zeta}{\pi}\right)dtd\zeta
\nonumber
\end{eqnarray}
and
\begin{eqnarray}
\lefteqn{
H_
{m, n}^{+, 2}(x)=
\frac{4iM^2}{\pi}\int\limits_{t=-\infty}^\infty
\int\limits_{\zeta=-T^\varepsilon}^{T^\varepsilon}t 
 e^{-t^2}V(m^2n, tM+T)\sin(x\cosh\zeta)e\left(\frac{tM\zeta}{\pi}\right)\nonumber}\\
&& \hspace{8.9cm}\times e\left(\frac{T\zeta}{\pi}\right)dtd\zeta.
\nonumber
\end{eqnarray}
In the following we only treat $H_{m, n}^{+, 1}(x)$ since $H_{m, n}^{+, 2}(x)$ is a lower order term which can be handled in a similar way.
 It is clear that
\begin{equation}
H_{m, n}^{+, 1}(x)=\frac{4iMT}{\pi}\int\limits_{\zeta=-T^\varepsilon}^{T^\varepsilon}\hat{k^*}
\left(-\frac{M\zeta}{\pi}\right)\sin(x\cosh\zeta)e\left(\frac{T\zeta}{\pi}\right)d\zeta
\end{equation}
which is equal to
$$
4iT\int\limits_{\zeta=-\frac{MT^\varepsilon}{\pi}}^{\frac{MT^\varepsilon}{\pi}}\hat{k^*}(\zeta)\sin\left(x\cosh\frac{\zeta\pi}{M}\right)e\left(-\frac{T\zeta}{M}\right)d\zeta
$$
by making a change of variable $-\frac{M\zeta}{\pi}\rightarrow\zeta,$ here 
\begin{equation}
k^*(t)=e^{-t^2}V(m^2n, tM+T)
\end{equation}
and
\begin{equation}
\hat{k^*}(\zeta)=\int\limits_{-\infty}^\infty k^*(t)e(-t\zeta)dt
\end{equation}
is its Fourier transform. Since $\hat{k^*}(\zeta)$ is a Schwartz class function, one can extend the integral in (4.10) to $(-\infty, \infty)$ with a negligible error term.
 Now let
\begin{equation}
W_{m, n}(x):=T\int\limits_{-\infty}^\infty \hat{k^*}(\zeta)\sin\left(x\cosh\frac{\zeta\pi}{M}\right)e\left(-\frac{T\zeta}{M}\right)d\zeta
\end{equation}
and
\begin{equation}
W_{m, n}^*(x):=T\int\limits_{-\infty}^\infty \hat{k^*}(\zeta)
e\left(-\frac{T\zeta}{M}-\frac{x}{2\pi}\cosh\frac{\zeta\pi}{M}\right)d\zeta,
\end{equation}
then
$$
W_{m, n}(x)=\frac{W_{m, n}^*(-x)-W_{m, n}^*(x)}{2i}
$$
and
$$
H_{m, n}^{+, 1}(x)=4iW_{m, n}(x)+O(T^{-A})
$$
with $A$ arbitrarily large. The contributaion to $W_{m, n}(x)$ coming from $|\zeta|\geqslant T^\varepsilon$ ($\varepsilon>0$ arbitrarily small but
 fixed) is negligible. So we need only consider $|\zeta|\leqslant T^\varepsilon.$  The phase $\phi$ in the exponential of $W_{m, n}^*(x)$  is
$$
\phi(\zeta)=-\frac{T\zeta}{M}-\frac{x}{2\pi}\cosh\frac{\zeta\pi}{M},
$$
so
$$
\phi^{'}(\zeta)=-\frac{T}{M}-\frac{x}{2M}\sinh\frac{\zeta\pi}{M}.
$$
Then if $|x|\leqslant T^{1-\varepsilon}M,$ $W_{m, n}^*(x)$ is negligible.  In the following we assume that
$$T^{1-\varepsilon}M\leqslant 
|x|\leqslant M^4.$$
In this case we need the asymptotic expansion of $W_{m, n}^*(x).$ One could quote Lemma 5.1 of [LYL]. For completeness,
we prefer to derive it here. But the methods are really based on [Sa] and [LYL]. Now
\begin{eqnarray}
\lefteqn{
W
_{m, n}^*(x)=T\int\limits_{-\infty}^\infty\hat{k^*}(\zeta)e\left(-\frac{T\zeta}{M}-\frac{x}{2\pi}-\frac{\pi x\zeta^2}{4M^2}
-\frac{\pi^3x\zeta^4}{48M^4}-\frac{\pi^5x\zeta^6}{1440M^6}\right)d\zeta}\nonumber\\
&&\hspace{6cm}+O\left(T\int
\limits_{-\infty}^\infty|\hat{k^*}(\zeta)|\frac{|\zeta|^8|x|}{M^8}d\zeta\right).
\nonumber
\end{eqnarray}
Expanding $e\left(\frac{-\pi^5x\zeta^6}{1440M^6}\right)$ into a Taylor series of order $1$, we have
\begin{equation}
W_{m, n}^*(x)=W_{m, n}^+(x)-\frac{2\pi^6ix}{1440M^6}W_{m, n}^-(x)+O\left(\frac{T|x|}{M^8}\right),
\end{equation}
where																																																																																										
$$
W_{m, n}^+(x)=Te\left(\frac{-x}{2\pi}\right)\int\limits_{-\infty}^\infty k_0^*(\zeta)e\left(-\frac{T\zeta}{M}-\frac{\pi x\zeta^2}{4M^2}\right)d\zeta
$$
with
$$
k_0^*(\zeta)=\hat{k^*}(\zeta)e\left(\frac{-\pi^3x\zeta^4}{48M^4}\right)
$$
and
$$
W_{m, n}^-(x)=Te\left(\frac{-x}{2\pi}\right)\int\limits_{-\infty}^\infty k^*
_
1(\zeta)e\left(-\frac{T\zeta}{M}-\frac{\pi x\zeta^2}{4M^2}\right)d\zeta
$$
with
$$
k_1^*(\zeta)=\zeta^6\hat{k^*}(\zeta)e\left(\frac{-\pi^3x\zeta^4}{48M^4}\right).
$$
Now by completing the square, we have
$$
W_{m, n}^+(x)=Te\left(\frac{-x}{2\pi}+\frac{T^2}{\pi x}\right)\int\limits_{-\infty}^\infty k_0^*(\zeta)e\left(-\frac{\pi x}{4M^2}\left(\zeta+\frac{2MT}{\pi x}\right)^2\right)d\zeta
$$
which is equal to ([GR], 3.691 1)
$$
(1+i)Te\left(\frac{-x}{2\pi}+\frac{T^2}{\pi x}\right)\int\limits_{-\infty}^\infty \hat{k_0^*}(\zeta)e\left(\frac{-2MT\zeta}{\pi x}\right)\frac{M}{\sqrt{\pi |x|}}e\left(\frac{M^2\zeta^2}{\pi x}\right) d\zeta
$$
by Parseval. Expanding $e\left(\frac{M^2\zeta^2}{\pi x}\right)$ in a Taylor series we have
\begin{eqnarray}
\lefteqn{
W_{m, n}^+(x)=(1+i)\frac{TM}{\sqrt{\pi |x|}}e\left(\frac{-x}{2\pi}+\frac{T^2}{\pi x}\right)
}\nonumber\\
&&\hspace{3cm}\times\sum\limits_{l=0}^{\infty}
\frac{1}{l!}\left(\frac{2iM^2}{x}\right)^l\int\limits_{-\infty}^\infty\zeta^{2l}\hat{k_0^*}(\zeta)e\left(\frac{-2MT\zeta}{\pi x}\right)d\zeta
 \nonumber\\
&&\hspace{0.7cm}=(1+i)\frac{TM}{\sqrt{\pi |x|}}e\left(\frac{-x}{2\pi}+\frac{T^2}{\pi x}\right)
\sum\limits_{l=0}^{\infty}\frac{(2i)^{-l}}{l!}\left(\frac{M^2}{\pi^2x}\right)^l{k_0^*}^{(2l)}\left(\frac{-2MT}{\pi x}\right)
.\nonumber
\end{eqnarray}
Since
$$
{k_0^*}^{(2l)}(t)=\sum\limits_{0\leqslant l_1\leqslant 2l}
\left(\!
\begin{matrix}
2l\\
l_1
\end{matrix}
\!\right)
\frac{d^{l_1}}{dt^{l_1}}e\left(\frac{-\pi^3xt^4}{48M^4}\right)\times 
\frac{d^{2l-l_1}}{dt^{2l-l_1}}\hat{k^*}(t)
$$
where 
$\begin{pmatrix}
n\\
r
\end{pmatrix}$ denotes the binomial coefficient and 
$$
\frac{d^{l_1}}{dt^{l_1}}e\left(\frac{-\pi^3xt^4}{48M^4}\right)\bigg|_{t=\frac{-2MT}{\pi x}}\ll 1,
$$ one can truncate the above series of $W_{m, n}(x)$ at order $L_1$ with a reminder $O\left(T\left(\frac{M}{\sqrt{x}}\right)^{2L_1+3}\right).$
Now expanding $e\left(\frac{-\pi^3xt^4}{48M^4}\right)$ in a power series and differentiating it termwisely, we have
\begin{eqnarray}
\lefteqn{
\frac{d^{l_1}}{dt^{l_1}}e\left(\frac{-\pi^3xt^4}{48M^4}\right)\bigg|_{t=\frac{-2MT}{\pi x}}\!\!
=\!\!\sum\limits_{4l_2\geqslant l_1}\frac
{(4l_2)!}{(4l_2-l_1)!l_2!}\left(\frac{2i\pi^4}{48}\right)^{l_2}
\left(\frac{-x}{M^4}\right)^{l_2}t^{4l_2-l_1}\bigg|_{t=\frac{-2MT}{\pi x}}}\nonumber\\
&&\hspace{3.2cm}=\!\!\!\!\sum\limits_{\frac{l_1}{4}\leqslant l_2\leqslant L_2}\!\!\!\!\frac
{(4l_2)!}{(4l_2-l_1)!l_2!}\!\!\left(\frac{i\pi^4}{24}\right)^{l_2}\!\!
\left(\frac{-x}{M^4}\right)^{l_2}\!\!\left(\frac{-2MT}{\pi x}\right)^{4l_2-l_1}\nonumber\\
&&\hspace{6.9cm}+O\left(\left(\frac{T^4}{|x|^3}\right)^{L_2+1}\left(\frac{|x|}{MT}\right)^{l_1}\right).
\nonumber
\end{eqnarray}
Combining the above, we have the following asymptotic expansion
\begin{eqnarray}
\lefteqn{
W_{m, n}^+(x)=\frac{TM}{\sqrt{|x|}}e\left(\frac{-x}{2\pi}+\frac{T^2}{\pi x}\right)\sum\limits_{l=0}^{L_1}
\sum\limits_{0\leqslant l_1\leqslant 2l}\sum\limits_{\frac{l_1}{4}\leqslant l_2\leqslant L_2}
c_{l, l_1, l_2}
}\nonumber\\
&&\hspace{4cm} \times \frac{M^{2l-l_1}T^{4l_2-l_1}}{x^{l+3l_2-l_1}}
{\hat{k^*}}^{(2l-l_1)}\left(\frac{-2MT}{\pi x}\right)\nonumber\\
&& \hspace{4cm}+ O\left(\frac{TM}{\sqrt{|x|}}\left(\frac{T^4}{|x|^3}\right)^{L_2+1}
 +T\left(\frac{M}{\sqrt{|x|}}\right)^{2L_1+3}\right),
\nonumber
\end{eqnarray}
here $c_{l, l_1, l_2}$ are constants depending only on $l, l_1$ and $l_2.$ 
$W_{m, n}^-(x)$ has similar asymptotic expansion. 
We end up with the following proposition (recall (4.15)):
\begin{proposition}{\rm  1)} For $|x|\leqslant T^{1-\varepsilon}M$ with $\varepsilon>0,$
$$
W_{m, n}^*(x)\ll T^{-A}
$$
where $A>0$ is arbitrarily large and the implied constant depends on $\varepsilon$ and $A$.\\
{\rm 2)} For $T^{1-\varepsilon}M\leqslant |x|\leqslant M^4, T^{\frac{3}{8}+\varepsilon}\leqslant M\leqslant
 T^{\frac{1}{2}}$ and $L_2, L_1\geqslant 1,$
\begin{eqnarray}
\lefteqn{\hspace{0.5cm}
W_{m, n}^*(x)=\frac{TM}{\sqrt{|x|}}e\left(\frac{-x}{2\pi}+\frac{T^2}{\pi x}\right)
\sum\limits_{l=0}^{L_1}\sum\limits_{0\leqslant l_1\leqslant 2l}
 \sum\limits_{\frac{l_1}{4}\leqslant l_2\leqslant L_2}
c_{l, l_1, l_2}\frac{M^{2l-l_1}T^{4l_2-l_1}}{x^{l+3l_2-l_1}}}\nonumber\\
&&\hspace{1.5cm}\times 
\left[{\hat{k^*}}^{(2l-l_1)}\left(\frac{-2MT}{\pi x}\right)
-\frac{2\pi^6ix}{1440M^6}(y^6\hat{k^*}(y))^{(2l-l_1)}\left(\frac{-2MT}{\pi x}\right)\right]\\
&&\hspace{3cm}+O\left(\frac{TM}{\sqrt{|x|}}\left(\frac{T^4}{|x|^3}\right)^{L_2+1}
+T\left(\frac{M}{\sqrt{|x|}}\right)^{2L_1+3}+\frac{T|x|}{M^8}\right)\nonumber
\end{eqnarray}
where $c_{l, l_1, l_2}$ are constants depending only on $l, l_1$ and $l_2,$ especially $c_{0, 0, 0}=\frac{1+i}{\sqrt{\pi}}.$
\end{proposition}
It follows from 1) in the above proposition, ${\cal R}_2^+$ is negligible.
The remaining part of this section is devoted to the estimation of ${\cal R}^+_3.$ Applying the asymptotic expansion (4.16)
of $W_{m, n}^*(x)$ and choosing $L_2$ and $L_1$ sufficiently large makes the contribution to ${\cal R}_3^+$ from the first two terms in the error term in (4.16) negligible. The contribution to
 ${\cal R}_3^+$ from the last term in the error term in (4.16) is
$$
O_{\varepsilon, f}\left(\frac{T^{1+\varepsilon}N}{M^8}\right)=O_{\varepsilon, f}\left(T^{1+\varepsilon}M\right)
$$
as expected, where we used the trivial bound for the Kloosterman sum and (2.7).  Since $|x|\geqslant T^{1-\varepsilon}M$, 
$$
\frac{M^{2l-l_1}T^{4l_2-l_1}}{x^{l-l_1+3l_2}}
\ll \left( \frac{M}{T^{1-\varepsilon}}\right)^l\left(\frac{T}{M^3}\right)^{l_2}T^{(3l_2-l_1)\varepsilon}\ll 1.
$$  From now on, we only take the leading term $l=0, l_1=0$ and $l_2=0$ in (4.16). The other terms are of an identical form and can be treated similarly.  We are led to estimate
                                                                                                                                                                            \begin{eqnarray}
\lefteqn{\;\;\;
\tilde{\cal R}^+_3:=\sqrt{2}i\pi^{-1}MTe\left(-\frac{1}{8}\right)
\sum\limits_{m\geqslant 1}\sum\limits_{n\geqslant 1}\frac{A(n, m)}{mn^{\frac{3}{4}}}g\left(\frac{m^2n}{N}\right)}\\
&&\hspace{1cm}\times\sum\limits_{c\leqslant C_2/m}c^{-\frac{1}{2}}S(n, 1; c)e\left(\frac{2\sqrt{n}}{c}-\frac{T^2c}{4\pi^2\sqrt{n}}\right)\hat{k^*}\left(\frac{MTc}{2\pi^2\sqrt{n}}\right)\nonumber.\nonumber
\end{eqnarray}
In the above, if we sum over $n$ trivially and applying Weil's bound for the Kloosterman sum
$$
S(n, 1; c)\ll_\varepsilon c^{\frac{1}{2}+\varepsilon},
$$
we have
$$
\tilde{\cal R}^+_3\ll MTC_2^{1+\varepsilon}N^{\frac{1}{4}}\ll T^{\frac{9}{4}+\varepsilon}.
$$
To save $T^{\frac{5}{4}}M^{-1},$ we have to sum over $n$ nontrivially by the Voronoi formula for $GL(3)$ (i.e, 
Proposition 2.1). 
Expanding the Kloosterman sum in (4.17) 
and applying
Proposition 2.1 with
$$
\psi(y)=y^{-\frac{3}{4}}g\left(\frac{m^2y}{N}\right)e\left(\frac{2\sqrt{y}}{c}-\frac{T^2c}{4\pi^2\sqrt{y}}\right)\hat{k^*}\left(\frac{MTc}{2\pi^2\sqrt{y}}\right),
$$
we have
\begin{eqnarray}
\lefteqn{
\sum\limits_{n\geqslant 1}A(n, m)e\left(\frac{n\bar{d}}{c}\right)\psi(n)}\nonumber\\
&&=\frac{c\pi^{-\frac{5}{2}}}{4i}\sum\limits_{n_1|cm}\sum\limits_{n_2>0}\frac{A(n_2, n_1)}{n_1n_2}
S(md, n_2; mcn_1^{-1})\Psi_{0, 1}^0\left(\frac{n_2n_1^2}{c^3m}\right)
\nonumber\\
&&\;\;+ \frac{c\pi^{-\frac{5}{2}}}{4i}\sum\limits_{n_1|cm}\sum\limits_{n_2>0}\frac{A(n_2, n_1)}{n_1n_2}
S(md, -n_2; mcn_1^{-1})\Psi_{0, 1}^1\left(\frac{n_2n_1^2}{c^3m}\right)
\nonumber
\end{eqnarray}
where $\Psi_{0, 1}^0(x)$ and $\Psi_{0, 1}^1(x)$ are defined below (2.5). As we expalined before Proposition 2.1, 
we only consider the first term involving $\Psi_0(x)$ on the right side of the above formula since all the other terms can be treated in a similar way.
Since $c\leqslant \frac{N^{\frac{1}{2}}}{T^{1-\varepsilon}Mm},$
$$
\frac{n_2n_1^2}{c^3m}\frac{N}{m^2}\gg T^{\frac{3}{2}},
$$
by Lemma 2.1 for $x=\frac{n_2n_1^2}{c^3m},$
\begin{eqnarray}
\lefteqn{\;\;\;\Psi_0(x)=2\pi^4xi\int\limits_0^\infty\psi(y)\frac{d_1\sin (6\pi x^{\frac{1}{3}}y^{\frac{1}{3}})}{(\pi^3xy)^{\frac{1}{3}}}dy +\text{lower order terms}}\\
&&\;\;\;\;\;\;=
\pi^3x^{\frac{2}{3}}d_1\int\limits_0^\infty e(u_1(y))a(y)dy-\pi^3x^{\frac{2}{3}}d_1\int\limits_0^\infty e(u_2(y))a(y)dy \nonumber\\
&&\hspace{6cm}\;\;\;+\text{lower order
 terms}\nonumber
\end{eqnarray}
where
$$
u_1(y)=\frac{2\sqrt{y}}{c}+3x^{\frac{1}{3}}y^{\frac{1}{3}},
$$
$$
u_2(y)=\frac{2\sqrt{y}}{c}-3x^{\frac{1}{3}}y^{\frac{1}{3}}
$$
and 
$$
a(y)=g\left(\frac{m^2y}{N}\right)
\hat{k^*}\left(\frac{MTc}{2\pi^2\sqrt{y}}\right)
e\left(-\frac{T^2c}{4\pi^2\sqrt{y}}\right)y^{-\frac{13}{12}}.
$$
Since $u_1^{'}
(y)\gg c^{-1}y^{-\frac{1}{2}}$ and $a^{'}
(y)\ll T^2cy^{-\frac{31}{12}},$ we have 
$$ u_1^{'}(y)
a^{'}
(y)^{-1}\gg M^2T^{-\varepsilon}\gg T^{\frac{3}{4}
-\varepsilon}.
$$
By partial integration many times, one shows that the contribution to (4.17) from the first integral in (4.18) is negligible.\\
Now we turn to the second integral in (4.18). Since
$$
u_2^{'}
(y)=\frac{1}{c}\sqrt{\frac{1}{y}}-x^{\frac{1}{3}}
y^
{-\frac{2}{3}},
$$
if
\begin{equation}
x\geqslant 2\frac{N^{\frac{1}{2}}}{mc^3} \;\;\;\;\;\;
\text {or} \;\;\;\;\;\;
x\leqslant \frac{2}{3}\frac{
N^{\frac{1}{2}}}
{
mc^3},
\end{equation}
then $$ u_2^{'}
(y)\gg \frac{1}{c}\sqrt{\frac{1}{y}}.
$$
As the argument above, under the condition (4.19), the contribution to (4.17) from the second integral in (4.18) is also negligible. So for stationary or small values
of $u_2^{'}
(y)$ 
we need only consider the case when
\begin{equation}
\frac{2}{3}\frac{
N^{\frac{1}{2}}}{
mc^3}\leqslant x\leqslant 2\frac{N^{\frac{1}{2}}}
{
mc^3}, \;\;\;\;\;
\text{i.e.,} 
\;\;\;\;\;
\frac{2}{3}\frac{
N^{\frac{1}{2}}}{
n_1^2}\leqslant
 n_2
\leqslant 2\frac{N^{\frac{1}{2}}}
{n_1^2}.
\end{equation}
Then 
$$
\int\limits_0^\infty e(u_2(y))a(y)dy=\int\limits_{\frac{x^2c^6}{4}}^{\frac{9}{2}x^2c^6}e(u_2(y))a(y)dy.
$$
There is a stationary phase point $y_0=x^2c^6$ such that $u_2^{'}(y_0)=0.$ Applying the stationary phase method ([Hu], p. 114), we have
\begin{eqnarray}
\lefteqn{
\int\limits_0^\infty e(u_2(y))a(y)dy}\\
&&=\frac{e\left(-xc^2+\frac{1}{8}\right)a(y_0)}{\sqrt{u_2^{''}(y_0)}}+O\left(c^{\frac{7}{2}}T^4N^{-\frac{11}{6}}m^{\frac{11}{3}}\right).
\nonumber
\end{eqnarray}
Due to 
\begin{equation}
\sum_{\substack{0\leqslant d<c\\ (d, c)=1}}e\left(\frac{
d}{c}\right)S(md, n_2; mcn_1^{-1})=\sum
_{\substack{u
(
\text{mod}\;mcn_1^{-1})\\u\bar{u}\equiv 1 (
\text{mod}\;mcn_1^{-1})}}
S(0, 1
+un_1; c)e\left(\frac{n_2\bar{u}}{mcn_1^{-1}}\right)
\end{equation}
where
$$
S(0, a; c)=\sum_{\substack{v
(\text{mod}\; c)\\(v
, c)=1}}e\left(\frac{av
}{c}\right)
$$
is the Ramanujan sum which is bounded by $(a, c),$ we deduce that (4.22) is bounded by $mc^{1+\varepsilon}$ with $\varepsilon>0.$
Therefore, the contribution to (4.17) from the error term in (4.21) is bounded by 
\begin{eqnarray}
\lefteqn{
MT\sum\limits_{m\geqslant 1}m^{-1}
\sum\limits_{c\leqslant C_2/m}c^{\frac{1}{2}}
\sum\limits_{n_1|cm}
\sum\limits_{\frac{2N^{\frac{1}{2}}}{3n_1^2}\leqslant n_2\leqslant 2\frac{N^{\frac{1}{2}}}{n_1^2}}
\frac{|A(n_1, n_2)|}{n_1n_2}
}\\
&&\hspace{2cm}
\times\left
(\frac{n_2n_1^2}{c^3m}\right)^{\frac{2}{3}}(mc)^{1+\varepsilon}c^
{\frac{7}{2}}T^4N^{-\frac{11}{6}}m^{\frac{11}{3}}\nonumber\\
&&\hspace{5cm}
\ll M^{-3}T^{1+\varepsilon}N^{\frac{1}{2}}
\ll T^{1+\varepsilon}M
\nonumber\end{eqnarray}
because $M\geqslant T^{\frac{3}{8}}.$ 
We conclude from (4.17), (4.21), (4.22) and (4.23) that
\begin{eqnarray}
\lefteqn{\;\;\;\;\;\;
\tilde{\cal R}^+_3=\pi^{-1}MT\sum\limits_{m\geqslant 1}m^{-1}\sum\limits_{c\leqslant C_2/m}c^{-1}\sum\limits_{n_1|cm}n_1^{-1}\sum\limits_{n_2>0}A(n_1, n_2)}\\
&&\;\;\;\;\;\;\;\;\;\times 
\sum\limits_
{\substack{0\leqslant u<m
cn_1^{-1}
\\u\bar{u}\equiv 1(\text{mod} \;
mcn_1^{-1}
)
}}S(0, 1+un_1; c)e\left(\frac{n_2\bar{u}}{mcn_1^{-1}}\right)e\left(\frac{-n_2n_1^2}{cm}\right)b(n_2)\nonumber\\
&&\hspace{8
cm}\;\;\;\;\;\;\;
+O(T^{1+\varepsilon}M)
\nonumber
\end{eqnarray}
where
$$
b(y)=y^{-1}g\left(\frac{y^2n_1^4}{N}\right)\hat{k^*}\left(\frac{MTcm}{2\pi^2yn_1^2}\right)e\left(
\frac{-T^2cm}{4\pi^2yn_1^2}\right).
$$
If we sum over $n_2$ trivially, we have
$$
{\cal R}_3^{+, *}\ll MT^{1+\varepsilon
}C_2\ll MT^{1+\varepsilon}\frac{T^{\frac{1}{2}}}{M}.
$$
In order
 to save $T^{\frac{1}{2}
}M
^{-1},$ we have to sum over $n_2$ nontrivially using the Voronoi formula for $GL(3)$ the second time. 
Invoking Proposition 2.1, one has 
\begin{eqnarray}
\lefteqn{
\sum\limits_{n_2
\geqslant 1}A(n
_1,n_2
)e\left(\frac{n_2(\bar{u}-n_1)}{mcn_1^{-1}}
\right)b(n_2)
}\\
&&=\frac{c^{'}
\pi^{-\frac{5}{2}}}{4i}\sum\limits_{l
_1|c^{'}n_1
}\sum\limits_{l
_2>0}\frac{A(l
_2, l
_1)}{l
_1l
_2}S(n_1\bar{u
^{'}}, l_2; n_1c^{'}l_1^{-1})
B
_{0, 1}^0\left(\frac{l
_2l
_1^2}{{
c^{'}}
^3
n_1
}\right)\nonumber\\
&&\;\;+\frac{c^{'}
\pi^{-\frac{5}{2}}}{4i}\sum\limits_{l
_1|c^{'}n_1
}\sum\limits_{l
_2>0}\frac{A(l
_2, l
_1)}{l
_1l
_2}S(n_1\bar{u
^{'}}, -l_2; n_1c^{'}l_1^{-1})
B
_{0, 1}^1\left(\frac{l
_2l
_1^2}{{
c^{'}}
^3
n_1
}\right)
\nonumber
\end{eqnarray}
where
$$
\frac{\bar{u}-n_1}{mcn_1^{-1}}:=\frac{u^{'}}{c^{'}}
$$ 
with $(u^{'}, c^{'})=1, c^{'}|mcn_1^{-1}$ and $B_{0, 1}^0(x)$ and $B_{0, 1}^1(x)$ are defined below (2.5). 
As before, we only consider the first term involving $B_{0, 1}^0(x)$  in (4.25)  since all the other terms can be treated in a similar way. Since
$$
\frac{l_2l_1^2}{{
c^{'}}
^3n_1}\frac{\sqrt{N}}{n_1^2}\gg
T^{1-\varepsilon}M,
$$
by Lemma 2.1 for $x=\frac{l_2l_1^2}{{
c^{'}}
^3n_1},$
\begin{eqnarray}
\lefteqn{B
_0(x)=2\pi^4xi\int\limits_0^\infty b
(y)\frac{d_1\sin (6\pi x^{\frac{1}{3}}y^{\frac{1}{3}})}{(\pi^3xy)^{\frac{1}{3}}}dy +\text{lower order terms}}\\
&&\;\;\;=
\pi^3x^{\frac{2}{3}}d_1\int\limits_0^\infty e(v
_1(y))q
(y)dy-\pi^3x^{\frac{2}{3}}d_1\int\limits_0^\infty e(v
_2(y))q
(y)dy \nonumber\\
&&\hspace{6cm}\;+\text{lower order
 terms}\nonumber
\end{eqnarray}
where
\begin{equation}
v_1(y)=3x^{\frac{1}{3}}y^{\frac{1}{3}}
-\frac{T^2cm}{4\pi^2yn_1^2},
\end{equation}
\begin{equation}
v_2(y)=-3x^{\frac{1}{3}}y^{\frac{1}{3}}
-\frac{T^2cm}{4\pi^2yn_1^2},
\end{equation}
and 
\begin{equation}
q(y)=y^{-\frac{4}{3}}g\left(\frac{y^2n_1^4}{N}\right)\hat{k^*}\left(\frac{MTcm}{2\pi^2yn_1^2}\right).
\end{equation}
Since
$$
v
_1^{'}(y)=x^{\frac{1}{3}}
y^{-\frac{2}{3}}+\frac{T^2cm}{4\pi^2y^2n_1^2}\gg\frac{T^2cm}{y^2n_1^2}
$$
and $q^{'}(y)\ll y^{-\frac{7}{3}}T^\varepsilon,$ we have
$$
v_1^{'}(y)q^{'}(y)^{-1}\gg y^{\frac{4}{3}}\frac{T^{2-\varepsilon}cm}{\sqrt{N}}\gg T^{\frac{5}{3}},
$$
by partial integration many times, one shows that the contribution to (4.24) from the first integral in (4.26) is negligible.\\
Now we turn to $v_2(y)$ defined by (4.28). Since
$$
v_2^{'}(y)=-x^{\frac{1}{3}}y^{-\frac{2}{3}}+\frac{T^2cm}{4\pi^2y^2n_1^2},
$$
if
\begin{equation}
x\geqslant \frac{T^6c^3m^3n_1^2}{10\pi^6N^2}
 \;\;\;\;\;\;
\text {or} \;\;\;\;\;\;
x\leqslant \frac{T^6c^3m^3n_1^2}{1000 \pi^6 N^2},
\end{equation}
one has
$$
|v_2^{'}(y)|\gg \frac{T^2cm}{y^2n_1^2}.
$$
As the arguments above, one shows that under the condition (4.30), the contribution to (4.24) from the second integral in (4.26) is negligible. For the remaining case
\begin{equation}
\frac{T^6c^3m^3n_1^2}{1000
 \pi^6 N^2}\leqslant x\leqslant \frac{T^6c^3m^3n_1^2}{10\pi^6N^2}, 
 \;\;\;\;\;
\text{i.e.,} 
\;\;\;\;\;\frac{L_2}{1000}\leqslant l_2\leqslant \frac{L_2}{10}
\end{equation}
with
$$
L_2=\frac{T^6c^3m^3n_1^3{c^{'}}^3}{\pi^6N^2l_1^2},
$$
we have
$$
|v_2^{''}(y)|\gg\frac{T^2cm}{y^3n_1^2}\gg T^2cmN^{-\frac{3}{2}}n_1^4.
$$
Therefore, by the second derivative test ([Hu], p. 88), one derives that
\begin{eqnarray}
\lefteqn{
B_0(x)\ll x^{\frac{2}{3}}\left(T^2cmN^{-\frac{3}{2}}n_1^4\right)^{-\frac{1}{2}}
\left(\frac{\sqrt{N}}{n_1^2}\right)^{-\frac{4}{3}}
T^\varepsilon}\\
&&\;\;\;
\ll T^{3+\varepsilon} c^{\frac{3}{2}}N^{-\frac{5}{4}}n_1^2m^{\frac{3}{2}}.\nonumber
\end{eqnarray}
Combining (4.24), (4.25), (4.32) and invoking the trivial bound for the Kloosterman sum one concludes that
\begin{eqnarray}
\lefteqn{
\tilde{\cal R}^+_3\ll MT\sum\limits_{m\geqslant 1}m^{-1}
\sum\limits_{c\leqslant C_2/m}c^{-1}\sum\limits_{n_1|cm}n_1^{-1}}\nonumber\\
&&\;\;\;\;\;\;
\times\sum
\limits_{u(\text{mod}
\;mcn_1^{-1})}(1+un_1, c) c^{'}
\sum\limits_{l_1|c^{'}n_1}\sum\limits_{\frac{L_2}{100}\leqslant l_2\leqslant \frac{L_2}{10}}
\frac{|A(l_1, l_2)|}{l_1l_2}\nonumber\\
&&\hspace{3cm}
\times n_1c^{'}l_1^{-1}T^{3+\varepsilon}c^{\frac{3}{2}}N^{-\frac{5}{4}}n_1^2m^{\frac{3}{2}}+O(MT^{1+\varepsilon})\nonumber\\
&&\;\;
\ll NT^{-\frac{1}{2}}M^{-\frac{7}{2}}+O(MT^{1+\varepsilon})\ll MT^{1+\varepsilon}\nonumber
\end{eqnarray}
since $M\geqslant T^{\frac{3}{8}}
.$ This finishes the estimation of ${\cal R}^+.$
\section{The terms related to the $K-$Bessel function}
\setcounter{equation}{0} 
This section is devoted to the estimation of ${\cal R}^-$ which is defined by (2.23). We split ${\cal R}^-$ into two parts 
${\cal R}_1^-$ and 
${\cal R}_2^-$ with
\begin{equation}
{\cal R}_1^-
= \sum\limits_{m\geqslant 1}\sum\limits_{n\geqslant 1}\frac{A(n, m)}{(m^2n)^{\frac{1}{2}}}g\left(\frac{m^2n}{N}\right)
\sum\limits_{c\geqslant C/m
}c^{-1}S(n, 1; c)H_{m, n}^-\left(\frac{4\pi\sqrt{n}}{c}\right)
\end{equation} and 
\begin{equation}
{\cal R}_2^-
= \sum\limits_{m\geqslant 1}\sum\limits_{n\geqslant 1}\frac{A(n, m)}{(m^2n)^{\frac{1}{2}}}g\left(\frac{m^2n}{N}\right)
\sum\limits_{c\leqslant C/m
}c^{-1}S(n, 1; c)H_{m, n}^-\left(\frac{4\pi\sqrt{n}}{c}\right),
\end{equation}
here
\begin{equation}
C=\sqrt{N}+T.
\end{equation}
First we will estimate (5.1). By (2.24) and the following formula ([Wa], p. 78)
$$
K_\nu(z)=\frac{1}{2}\pi\frac{I_{-\nu}(z)-I_\nu(z)}{\sin\nu\pi}
$$
where $I_\nu(z)$ is the $I-$Bessel function, we have
\begin{eqnarray}
\lefteqn{
H_{m, n}^-(x)=2\int\limits_{-\infty}^\infty\frac{I_{-2it}
(x)-I_{2it}(x)}{\sin 2it\pi}\sinh(\pi t)k(t)V(m^2n, t)tdt\nonumber}\\
&&\;\;\;\;\;\;\;
=-4\int\limits_{-\infty}^\infty\frac{I_{2it}(x)}{\sin 2it\pi}\sinh(\pi t) k(t)V(m^2n, t)tdt.\nonumber
\end{eqnarray}
Moving the line of integration to $\Im t=-\sigma=-100,$ $H_{m, n}^-(x)$ becomes
\begin{eqnarray}
\lefteqn{
-4\int\limits_{-\infty}^\infty \left[\sin\pi(2\sigma+2iy)\right]^{-1}I_{2\sigma+2iy}(x)\sinh\pi(-\sigma i+y)}\\
&\hspace{3cm}
\times k(-\sigma i+y)V(m^2n, -\sigma i+y)(-\sigma i+y)dy.\nonumber\end{eqnarray}
By the following formula ([GR], 8.431 3)
$$
I_\nu(x)=\frac{\left(\frac{x}{2}\right)^\nu}{\Gamma(\nu+\frac{1}{2})\Gamma(\frac{1}{2})}\int\limits_0^\pi e^{x\cos\theta}\sin^{2\nu}\theta d\theta
$$ for $\Re \nu>-\frac{1}{2},$ one derives that
\begin{equation}
I_{2\sigma+2iy}(x)\ll_\sigma x
^{2\sigma}|y|^{-2\sigma}e^{\pi y}e^x.
\end{equation}

Combining (5.4), (5.5) and (4.7), we have
\begin{equation}
H_{m, n}^-(x)\ll x^{2\sigma}e^x(m^2n)^{-\sigma}T^{\sigma+1+\varepsilon}M.
\end{equation}
By (2.7), (5.5) and the trivial bound for the Kloosterman sum, one obtains that
\begin{eqnarray}
\lefteqn{\hspace{-5cm}\;\;\;\;\;
{\cal R}_1^-\ll\sum\limits_{m\geqslant
 1}\sum\limits_{n\geqslant 1}
\frac{|A(n, m)|}{(m^2n)^{\frac{1}{2}}}
g\left(\frac{m^2n}{N}\right)\sum\limits_{c\geqslant C/m}
\left(\frac{\sqrt{n}}
{c}\right)^{2\sigma}
T^{\sigma+1+\varepsilon}
(m^2n)^{-\sigma}M}\nonumber\\
&\hspace{-5cm}
\ll N^{\frac{1}{2}}T^{2-\sigma+\varepsilon}M\ll 1.\nonumber
\end{eqnarray}
It remains to estimate ${\cal R}_2^-.$ By the following integral representation of the $K-$Bessel function (see [GR], 8.432 4)
$$
K_{2it}(x)=\frac{1}{2}\cosh ^{-1}t\pi \int\limits_{-\infty}^\infty \cos(x\sinh \zeta)e\left(-\frac{t\zeta}{\pi}\right)d\zeta
$$
and partial integration in $\zeta$ once, we have
\begin{eqnarray}
\lefteqn{
H^-_{m, n}(x)=\frac{4}{\pi}\int\limits_0^\infty\int\limits_{|\zeta|\leqslant T^\varepsilon}
\tanh \pi te^{-\frac{(t-T)^2}{M^2}}V(m^2n, t) t\cos (x\sinh\zeta)}\nonumber\\
&&\hspace{6cm}\times e\left(-\frac{t\zeta}{\pi}\right)d\zeta dt+O(T^{-A})\nonumber
\end{eqnarray}
where $A$ is arbitrarily large. By making change of a variable $\frac{t-T}{M}\rightarrow t,$
\begin{eqnarray}
\lefteqn{
H_{m, n}^-(x)=\frac{4M}{\pi}\int\limits_{-\frac{T}{M}}^\infty\int\limits_{|\zeta|\leqslant T^\varepsilon}\tanh \pi (tM+T)e^{-t^2}V(m^2n, tM+T)}\nonumber\\
&&\hspace{2.5cm} \times (tM+T)\cos (x\sinh \zeta)e\left(-\frac{tM\zeta}{\pi}-\frac{T\zeta}{\pi}\right)dtd\zeta+O(T^{-A}).
\nonumber
\end{eqnarray}
Following  the derivation of Proposition 4.1, by extending the $t$ integral to $(-\infty, \infty)$ with a negligible 
error term, we have
$$
H_{m, n}^-(x)=H_{m, n}^{-, 1}(x)+H_{m, n}^{-, 2}(x) + O(T^{-A}),
$$
where
\begin{eqnarray}
\lefteqn{
H_
{m, n}^{-, 1}(x)=
\frac{4MT}{\pi}\int\limits_
{t=-\infty}^\infty
\int\limits_{|\zeta|\leqslant T^\varepsilon} 
 e^{-t^2}V(m^2n, tM+T)\cos(x\sinh\zeta)}\nonumber\\
&&\hspace{6cm}\times e\left(-\frac{(tM+T)\zeta}{\pi}\right)dtd\zeta\nonumber
\end{eqnarray}
and
\begin{eqnarray}
\lefteqn{
H_
{m, n}^{-, 2}(x)=
\frac{4M^2}{\pi}\int\limits_{t=-\infty}^\infty
\int\limits_{|\zeta|\leqslant T^\varepsilon}t 
 e^{-t^2}V(m^2n, tM+T)\cos(x\sinh\zeta)}\nonumber\\
&&\hspace{6cm}\times e\left(-\frac{(tM+T)\zeta}{\pi}\right)dtd\zeta.\nonumber
\end{eqnarray}
In the following we only treat $H_{m, n}^{-, 1}(x).$ $H_{m, n}^{-, 2}(x)$ is a lower order term which can be
 handled in a similar way.
 It is clear that
$$
H_{m, n}^{-, 1}(x)=\frac{4MT}{\pi}\int\limits_{|\zeta|\leqslant T^\varepsilon}\hat{k^*}\left(\frac{M\zeta}{\pi}\right)
\cos(x\sinh\zeta)e\left(-\frac{T\zeta}{\pi}\right)d\zeta
$$
which is equal to
\begin{equation}
4T\int\limits_{|\zeta|\leqslant \pi^{-1}MT^\varepsilon}\hat{k^*}(\zeta)\cos\left(x\sinh\frac{\zeta\pi}{M}\right)e\left(-\frac{T\zeta}{M}\right)d\zeta
\end{equation}
by making a change of variable $\frac{M\zeta}{\pi}\rightarrow\zeta.$  Since $\hat{k^*}(\zeta)$ is a Schwartz class function, one can extend the above integral to $(-\infty, \infty)$ with a negligible error term. Now let
$$
Y_{m, n}(x):=T\int\limits_{-\infty}^\infty \hat{k^*}(\zeta)\cos\left(x\sinh\frac{\zeta\pi}{M}\right)e\left(-\frac{T\zeta}{M}\right)d\zeta
$$
and 
\begin{equation}
Y_{m, n}^*(x):=T\int\limits_{-\infty}^\infty \hat{k^*}(\zeta)
e\left(-\frac{T\zeta}{M}+\frac{x}{2\pi}\sinh\frac{\zeta\pi}{M}\right)d\zeta,
\end{equation}
then
\begin{equation}
Y_{m, n}(x)=\frac{Y_{m, n}^*(x)+Y_{m, n}^*(-x)}{2}
\end{equation}
and
$$
H_{m, n}^{-, 1}(x)=4Y_{m, n}(x)+O(T^{-A})
$$
with $A$ arbitrarily large.  Let
$$
\Omega(\zeta)=\frac{x\sinh\frac{\zeta\pi}{M}}{2\pi}-\frac{T\zeta}{M},
$$
then
$$
\Omega^{'}(\zeta)=\frac{x\cosh\frac{\zeta\pi}{M}}{2M}-\frac{T}{M}.
$$
Then if 
$$
|x|\leqslant \frac{1}{100}T\;\;\;\text{or}\;\;\; |x|\geqslant 100T
$$
then
$$
\Omega^{'}(\zeta)\gg \frac{T}{M}\gg T^\varepsilon,
$$
hence by partial integrations,
$$
Y_{m, n}^*(x)\ll T^{-A}
$$
with $A>0$ arbitrarily large. We are left with the case when
$$
\frac{1}{100}T\leqslant x\leqslant 100T,
$$
then
$$
\frac{x}{M^3}\ll T^{-\frac{1}{8}}
$$
(recall $M\geqslant T^{\frac{3}{8}}$). Now
\begin{eqnarray}
\lefteqn{
Y_{m, n}^*(x)
=T\int\limits_{-\infty}^\infty\hat{k^*}(\zeta )
e\left(-\frac{T\zeta}{M}+\frac{x\zeta}{2M}+\frac{\pi^2x\zeta^3}{12M^3}
+\frac{\pi^4x\zeta^5}{240M^5}\right)d\zeta}
\nonumber\\
&&\hspace{5cm}+O\left(T\int\limits_{-\infty}^\infty |\hat{k^*}(\zeta)|\frac{|\zeta|^7|x|}{M^7}d\zeta\right).\nonumber
\end{eqnarray}
Expanding $e\left(\frac{\pi^2x\zeta^3}{12M^3}+\frac{\pi^4x\zeta^5}{240M^5}\right)$ into a Taylor series of order 
$L_2,$ we have
\begin{eqnarray}
\lefteqn{
Y_{m, n}^*(x)=T\int\limits_{-\infty}^\infty \hat{k^*}(\zeta)e\left(-\frac{(2T-x)\zeta}{2M}\right)d\zeta}\nonumber\\
&&\hspace{2cm}\times\sum\limits_{l=0}^{L_2}
\sum\limits_{j=0}^l d_{j, l}\left(\frac{x\zeta^3}{M^3}\right)^j\left(\frac{x\zeta^5}{M^5}\right)^{l-j}d\zeta\nonumber\\
&&\hspace{6cm}+O\left(\frac{T|x|^{L_2+1}}{M^{3L_2+3}}+\frac{T|x|}{M^7}\right)\nonumber
\end{eqnarray}
where $d_{j, l}$ are constants coming from the Taylor expansion and especially \\$d_{0, 0}=1.$ Clearly
\begin{eqnarray}
\lefteqn{
Y_{m, n}^*(x)=T\sum\limits_{l=0}^{L_2}\sum\limits_{j=0}^l d_{j, l}
\frac{x^l}{M^{5l-2j}}{k^*}^{(5l-2j)}\left(\frac{x-2T}{2M}\right)(2\pi i)^{-5l+2j}}\nonumber\\
&&\hspace{6cm}+O\left(\frac{T|x|^{L_2+1}}{M^{3L_2+3}}+\frac{T|x|}{M^7}\right).
\nonumber
\end{eqnarray}
We end up with the following proposition
\begin{proposition}
{\rm 1) }For $|x|\geqslant 100 T$ or $x\leqslant \frac{1}{100}T,$
$$
Y_{m, n}^*(x)\ll T^{-A}
$$
where $A>0$ is arbitrarily large and the implied constant depends only on $A$.\\
{\rm 2)} For $\frac{1}{100}T\leqslant |x|\leqslant 100T$, $T^{\frac{3}{8}+\varepsilon}\leqslant M\leqslant T^{\frac{1}{2}}$   and $L_2\geqslant 1,$
\begin{eqnarray}
\lefteqn{
Y_{m, n}^*(x)=T\sum\limits_{l=0}^{L_2}\sum\limits_{j=0}^l b_{j, l}
\frac{x^l}{M^{5l-2j}}{k^*}^{(5l-2j)}\left(\frac{x-2T}{2M}\right)}\nonumber\\
&&\hspace{5cm}+O\left(\frac{T|x|^{L_2+1}}{M^{3L_2+3}}+\frac{T|x|}{M^7}\right),
\nonumber
\end{eqnarray}
where $b_{j, l}$ are constants depending only on $j$ and $l$, especially $b_{0, 0}=1.$
\end{proposition}
The contribution to ${\cal R}^-_2$ from the error term $O\left(\frac{T|x|}{M^7}\right)$ in the above proposition is $O(T^{1+\varepsilon}M)$ by (2.7) and the trivial bound for the Kloosterman sum. 
We always take $L_2$ sufficiently large such that the first error term in Proposition 5.1 2) is negligible. From now on 
we only take the leading term $l=0$ since all the other lower order terms can be handled similarly. Let
\begin{eqnarray}\lefteqn{}\\
&&\!\!\!\!\!\!\!\!\!\!{\tilde{\cal{R}}}^-_2:=T\sum\limits_{m\geqslant 1}\sum\limits_{n\geqslant 1}\frac{A(n, m)}{(m^2n)^{\frac{1}{2}}}g\left(\frac{m^2n}{N}\right)\!\!
\!\!\!\!\sum\limits
_{\frac{\sqrt{N}}{100Tm}\leqslant c\leqslant \frac{100\sqrt{N}}{Tm}}\!\!\!\!
c^{-1}
 S(n, 1; c)k^*\left(\frac{\frac{4\pi\sqrt{n}}{c}-2T}{2M}\right).\nonumber
\end{eqnarray}
If we sum over $n$ trivially and apply Weil's bound for the Kloosterman sum, one derives that
$$
{\tilde{\cal{R}}}^-_2\ll T^{\frac{1}{2}}N^{\frac{3}{4}+\varepsilon}\ll T^{\frac{11}{4}+\varepsilon}.
$$
To save $T^{\frac{7}{4}}M^{-1},$ we have to sum over $n$ nontrivially using the Voronoi formula for $GL(3).$ Expanding the Kloosterman sum in (5.10),
 by Proposition 2.1, we have
\begin{eqnarray}
\lefteqn{
\sum\limits_{n
\geqslant 1}A(n
,m
)e\left(\frac{n\bar{a}}{c}
\right)r(n)
}\\
&&=\frac{c
\pi^{-\frac{5}{2}}}{4i}\sum\limits_{n
_1|cm
}\sum\limits_{n_2>0}\frac{A(n
_2, n
_1)}{n
_1n
_2}S(ma, n_2; mcn_1^{-1})
R_{0, 1}^0\left(\frac{n_2n_1^2}{c^3m}\right)
\nonumber\\
&&\;\;+ \frac{c\pi^{-\frac{5}{2}}}{4i}\sum\limits_{n_1|cm}\sum\limits_{n_2>0}\frac{A(n_2, n_1)}{n_1n_2}
S(ma, -n_2; mcn_1^{-1})R_{0, 1}^1\left(\frac{n_2n_1^2}{c^3m}\right)
\nonumber
\end{eqnarray}
where
$$
r(y)=g\left(\frac{m^2y}{N}\right)k^*\left(\frac{\frac{4\pi\sqrt{y}}{c}-2T}{2M}\right)y^{-\frac{1}{2}}
$$
and $R_{0, 1}^0(x)$  and $R_{0, 1}^1(x)$ are defined below (2.5). As before, in the following, we only consider $R_0(x)$ since $x^{-1}R_1(x)$ has similar asymptotic behavior as of $R_0(x)$. Since
$$
\frac{n_2n_1^2}{c
^3m}\frac{N }{m^2}\gg
\frac{T^3}{N^{\frac{1}{2}}}\gg T^{\frac{3}{2}-\varepsilon},
$$
by Lemma 2.1 for $x=\frac{n_2n_1^2}{c^3m},$
$$
R
_0(x)=2\pi^4xi\int\limits_0^\infty r
(y)\frac{d_1\sin (6\pi x^{\frac{1}{3}}y^{\frac{1}{3}})}{(\pi^3xy)^{\frac{1}{3}}}dy +\text{lower order terms}.
$$
If $n_2\gg\frac{N^{\frac{1}{2}}T^\varepsilon}{M^3n_1^2},$ then
$$
x^{\frac{1}{3}}y^{-\frac{2}{3}}[r^{'}(y)]^{-1}\gg T^{\varepsilon}.
$$
By partial integration many times, one shows that the contribution to ${\tilde{\cal{R}}}^-_2$ from such terms is negligible. Next we assume 
$$n_2\ll\frac{N^{\frac{1}{2}}T^\varepsilon}{M^3n_1^2}.$$ Since $k^*(y)\ll(1+|y|)^{-A}$ for any $A>0,$ $r(y)$ is negligible unless
$$
\left|\frac{\frac{2\pi\sqrt{y}}{c}-T}{M}\right|\leqslant T^\varepsilon
$$ which implies that
$$
\frac{1}{4\pi^2}(Tc-T^\varepsilon Mc)^2\leqslant y\leqslant \frac{1}{4\pi^2}(Tc+T^\varepsilon Mc)^2,
$$ then
\begin{equation}
R_0(x)\ll x^{\frac{2}{3}}\left(\frac{N}{m^2}\right)^{-\frac{5}{6}}T^{1+\varepsilon}Mc^2.
\end{equation}
Combining (5.10), (5.11), (4.22) and (5.12), we have 
\begin{eqnarray}
\lefteqn{
{\tilde{\cal{R}}}^-_2\ll T\sum\limits_{m\leqslant \sqrt{N}}\frac{1}{m}\sum\limits
_{\frac{\sqrt{N}}{100Tm}\leqslant c\leqslant \frac{100\sqrt{N}}{Tm}}\sum\limits_{n_1|cm}\sum\limits_{n_2\ll\frac{N^{\frac{1}{2}}T^\varepsilon}{M^3n_1^2}}}\nonumber\\
&&\hspace{1cm}\times\frac{|A(n_1, n_2)|}{n_1n_2}mc^{1+\varepsilon}\left(\frac{n_2n_1^2}{c^3m}\right)^{\frac{2}{3}}\left(\frac{N}{m^2}\right)^{-\frac{5}{6}}T^{1+\varepsilon}
Mc^2\nonumber\\
&&\ll N^{\frac{1}{2}}M^{-1}T^\varepsilon\ll T^{1+\varepsilon}M\nonumber
\end{eqnarray}
since $M\geqslant T^{\frac{3}{8}}.$\\
This finishes the estimation of ${\cal R}^-$ and hence the proof of the main theorem.
\renewcommand{\theequation}{A.\arabic{equation}}
 \setcounter{equation}{0}
  \renewcommand{\thetheorem}{A.\arabic{theorem}}\renewcommand{\theproposition}{A.\arabic{proposition}}\renewcommand{\thecorollary}{A.\arabic{corollary}}
  \section*{Appendix}  
In this appendix, we consider the subconvexity problem of $L(\frac{1}{2}, f\times h)$ where $f$ is a self dual Hecke-Maass form for $SL(3, \mathbb{Z})$ and $h$ runs through holomorphic Hecke cusp forms of weight $k\geqslant 2$ and congruent to $0 (\text{mod}\; 4)$ for $SL(2, \mathbb{Z})$. This analogous problem was suggested by Peter Sarnak and we would like to thank him here.\\
Let ${\cal B}_k(SL(2, \mathbb{Z}))$ denote an orthogonal basis of holomorphic Hecke cusp forms of weight $k\equiv 0(\text{mod}\; 4)$ for $SL(2, \mathbb{Z}),$ each $h$ in ${\cal B}_k(SL(2, \mathbb{Z}))$ is normalized to have the first Fourier coefficient $a_h(1)$ equal to $1.$ Set
$$
\lambda_h(n)=\frac{a_h(n)}{n^{\frac{k-1}{2}}}.
$$
By Deligne [De],
$$
|\lambda_h(n)|\leqslant \tau(n).
$$
For $f$ a self dual Heke-Maass form of type $(\nu, \nu)$ for $SL(3, \mathbb{Z})$ with the Fourier-Whittaker expansion (2.4) and $h\in {\cal B}_k(SL(2, \mathbb{Z})),$ we define the Rankin-Selberg $L-$function
$$
L(s, f\times h)=\sum\limits_{m=1}^\infty\sum\limits_{n=1}^\infty\frac{\lambda_h(n)A(n, m)}{(m^2n)^s}.$$
It is entire and satisfies the functional equation
\begin{equation}
\Lambda(s, f\times h)=\Lambda(1-s, f\times h)
\end{equation}
where
\begin{eqnarray}
\lefteqn{\!\!
\Lambda(s, f\times h)=\pi^{-3s}
\Gamma\left(\frac{s+\frac{k-1}{2}-\alpha}{2}\right)
\Gamma\left(\frac{s+\frac{k-1}{2}-\beta}{2}\right)\Gamma\left(\frac{s+\frac{k-1}{2}-\gamma}{2}\right)}\nonumber\\
&&\hspace{0.9cm}\times
\Gamma\left(\frac{s+\frac{k+1}{2}-\alpha}{2}\right)\Gamma\left(\frac{s+\frac{k+1}{2}-\beta}{2}\right)\Gamma\left(\frac{s+\frac{k+1}{2}-\gamma}{2}\right)L(s, f\times h)
\nonumber
\end{eqnarray}
and 
\begin{equation}
\alpha=-3\nu+1,\;\;\beta=0,\;\;\gamma=3\nu-1.
\end{equation}
The above functional equation can be obtained by examining the template arising from the case of the minimal parabolic Eisenstein series for $GL(3)$ twisted by a cusp form in ${\cal B}_k(SL(2, \mathbb{Z}))$ (see [Gol], p. 315). Note the sign of the above functional equation is $+1$ because we restrict $k$ to be congruent to $0 (\text{mod}\; 4)$
(see [IK] p. 131 and [Iw1] p. 121). This is important because we need the uniformity of the sign of the functional equations of $L(\frac{1}{2}, f\times h)$ when applying the Petersson formula. The main theorem in this appendix is
\begin{theorem}
Let $f$ be a fixed self dual Hecke-Maass form for $SL(3, \mathbb{Z})$, then for $\varepsilon>0, K$ large and
$K^{\frac{3}{8}+\varepsilon}\leqslant M\leqslant K^{\frac{1}{2}},$ 
we have
$$
\sum\limits_{2\leqslant k\equiv 0 (\text {mod}\; 4)} e^{-\frac{(k-K)^2}{M^2}}\sum\limits_{h\in {\cal B}_k(SL(2, \mathbb{Z}))}L\left(\frac{1}{2}, f\times h\right)
\ll_{\varepsilon, f}K^{1+\varepsilon}M.
$$
\end{theorem}
As we explained in the introduction, Lapid's theorem applies which means that $L(\frac{1}{2}\times h)\geqslant 0.$ Due to this important property, we have
\begin{corollary}
Under the same assumptions as in the above theorem,
$$
L\left(\frac{1}{2}, f\times h\right)\ll_{\varepsilon, f} k^{\frac{11}{8}+\varepsilon}.
$$
\end{corollary}
The corresponding convexity bound for $L(\frac{1}{2}, f\times h)$ is 
$k^{\frac{3}{2}+\varepsilon}$ with $\varepsilon>0,$ so the above bound breaks the convexity
bound.  The rest of the paper is devoted to the proof of Theorem A.1.
As in Lemma 2.2, we have the following approximate functional equation for $L(s, f\times h):$ 
\begin{equation}
L\left(\frac{1}{2}, f\times h\right)=2\sum\limits_{m\geqslant 1}\sum\limits_{n\geqslant 1}
\frac{\lambda_j(n)A(n, m)}{(m^2n)^{\frac{1}{2}}}
U(m^2n, k)
\end{equation}
where
$$
U(y, k)=\frac{1}{2\pi i}\int\limits_{(1000)}y^{-u}F(u)\frac{\gamma(\frac{1}{2}+u, k)}{\gamma(\frac{1}{2}, k)}
\frac{du}{u}
$$ and 
\begin{eqnarray}
\lefteqn{\;\;\;\;\;\;\;\;\;\;\;
\gamma(s, k)=\pi^{-3s}
\Gamma\left(\frac{s+\frac{k-1}{2}-\alpha}{2}\right)
\Gamma\left(\frac{s+\frac{k-1}{2}-\beta}{2}\right)
\Gamma\left(\frac{s+\frac{k-1}{2}-\gamma}{2}\right)}\nonumber\\
&&\hspace{2.5cm}\times
\Gamma\left(\frac{s+\frac{k+1}{2}-\alpha}{2}\right)
\Gamma\left(\frac{s+\frac{k+1}{2}-\beta}{2}\right)
\Gamma\left(\frac{s+\frac{k+1}{2}-\gamma}{2}\right).
\nonumber
\end{eqnarray}
We introduce the spectrally normalized first moment of the central values of $L-$functions
$$
{\cal A}:=\sum\limits_{2\leqslant k\equiv 0 (\text {mod}\; 4)} e^{-\frac{(k-K)^2}{M^2}}\sum\limits_{h\in {\cal B}_k(SL(2, \mathbb{Z}))}
 \frac{KL(\frac{1}{2}, f\times h)}{(k-1)L(1, \text{sym}^2h)}.
$$
The weights $L^{-1}(1, \text{sym}^2h)$ are needed in the Petersson formula and they are harmless since it is known ([Iw], [HL]) that
$$
k^{-\varepsilon}\ll L(1, \text{sym}^2h)\ll k^\varepsilon
$$
for any $\varepsilon>0.$ Applying (A.1) to ${\cal A},$ it is enough to show 
\begin{eqnarray}
\lefteqn{
\sum\limits_{2\leqslant k\equiv 0 (\text {mod}\; 4)} e^{-\frac{(k-K)^2}{M^2}}\frac{K}{k-1}\sum\limits_{m\geqslant 1}
\sum\limits_{n\geqslant 1}\frac{A(n, m)}{(m^2n)^{\frac{1}{2}}}
U(m^2n, k)}
\nonumber\\
&&\hspace{5cm}\times g\left(\frac{m^2n}{N}\right){\cal F}_k\ll K^{1+\varepsilon}M,
\end{eqnarray}
here $g$ is a fixed smooth function of compact support on $[1, 2]$, $1\leqslant N\ll_\varepsilon K^{3+\varepsilon}$ and
$$
{\cal{F}}_k=\sum\limits_{h\in {\cal B}_k(SL(2, \mathbb{Z}))}
\frac{\lambda_h(n)}{L(1, \text{sym}^2h)}.
$$
By Petersson's formula (see [ILS], p. 111, for example), 
\begin{equation}
{\cal F}_k=\frac{k-1}{2\pi^2}\left[\delta(n, 1)+2\pi\sum\limits_{c\geqslant 1}c^{-1}S(n, 1; c) 
J_{k-1}\left(\frac{4\pi\sqrt{n}}{c}\right)\right].
\end{equation}
We then write the left side of (A.4) as 
$${\cal D}_w+{\cal {ND}}_w,
$$
where
\begin{equation}
{\cal D}_w=\sum\limits_{2\leqslant k\equiv 0 (\text {mod}\; 4)} \frac{K}{2\pi^2}e^{-\frac{(k-K)^2}{M^2}}\sum\limits_{m\geqslant 1}
\frac{A(1, m)}{m}U(m^2, k)g\left(\frac{m^2}{N}\right)
\end{equation}
and 
\begin{eqnarray}
\lefteqn{\;\;\;\;\;\;\;\;\;
{\cal{ND}}_w=\sum\limits_{2\leqslant k\equiv 0 (\text {mod}\; 4)} \frac{K}{\pi}e^{-\frac{(k-K)^2}{M^2}}\sum\limits_{m\geqslant 1}
\sum\limits_{n\geqslant 1}\frac{A(n, m)}{(m^2n)^{\frac{1}{2}}}U(m^2n, k)
}\\
&& \hspace{4cm}\times g\left(\frac{m^2n}{N}\right)\sum\limits_{c\geqslant 1}c^{-1}
S(n, 1; c) J_{k-1}\left(\frac{4\pi\sqrt{n}}{c}\right).\nonumber
\end{eqnarray}
From (3.2),
$$
{\cal D}_w\ll K^{1+\varepsilon} M,
$$
which is consistent with the desired bound in (A.2).\\
To estimate ${\cal {ND}}_w,$ we begin by executing the $k$-sum by Poissson summation as in [Iw1] (p. 86) and [Sa] (p. 430). 
Applying the following integral representation [GR] of the $J$-Bessel function
$$
J_l(x)=\int\limits_{-\frac{1}{2}}^{\frac{1}{2}}e(lt)e^{-ix\sin 2\pi t} dt
$$ and the Poisson summation in $k$ yields 
\begin{equation}
K\sum\limits_{2\leqslant k\equiv 0 (\text {mod}\; 4)} u(k-1)J_{k-1}(x)=-\frac{1}{2}V_1(x)+\frac{i}{2}V_2(x)
\end{equation}
where
$u(x)=e^{-\frac{(x+1-K)^2}{M^2}}U(m^2n, x+1),$
\begin{equation}
V_1(x)=K\int\limits_{-\infty}^\infty\hat{u}(t)\sin (x\cos 2\pi t)dt,
\end{equation}
and 
\begin{equation}
V_2(x)=K\int\limits_{-\infty}^\infty\hat{u}(t)\sin (x\sin 2\pi t)dt,
\end{equation}
with $\hat{u}(t)$ be the Fourier transform of $u(x)$ as defined in (4.12). Since
$$
\hat{u}(t)=Me(-(K-1)t)\hat{u}_0(Mt)
$$
with
\begin{equation}
u_0(x)=e^{-x^2}U(m^2n, xM+K),
\end{equation}we have
\begin{equation}
V_1(x)=K\int\limits_{-\infty}^\infty\hat{u}_0(t)e\left(\frac{-(K-1)t}{M}\right)\sin \left(x\cos\frac{ 2\pi t}{M}
\right)dt,
\end{equation}
and 
\begin{equation}
V_2(x)=K\int\limits_{-\infty}^\infty\hat{u}_0(t)e\left(\frac{-(K-1)t}{M}\right)\sin 
\left(x\sin\frac{ 2\pi t}{M}\right)dt.
\end{equation}
We will first estimate the contribution to (A.7) from $V_1(x).$ Set 
$$
V_1^*(x)=K\int\limits_{-\infty}^\infty\hat{u}_0(t)e\left(\frac{-(K-1)t}{M}-\frac{x}{2\pi}\cos\frac{2\pi t}{M}\right)dt,
$$
then
$$
V_1(x)=\frac{V_1^*(-x)-V_1^*(x)}{2i}.\nonumber
$$
One can see that $V_1^*(x)$ and $W_{m, n}^*(x)$ (see (4.14)) have similar integral representation. 
Following the derivation of Proposition 4.1, it is straightforward  to derive the following:
\begin{proposition} {\rm 1)} For $|x|\leqslant K^{1-\varepsilon}M$ with $\varepsilon>0,$
$$
V_1^*(x)\ll K^{-A}
$$
where $A>0$ is arbitrarily large and the implied constant depends on $\varepsilon$ and $A$.\\
{\rm 2)} For $K^{1-\varepsilon}M\leqslant |x|\leqslant M^4, K^{\frac{3}{8}+\varepsilon}\leqslant M\leqslant 
K^{\frac{1}{2}}$ and $L_2, L_1\geqslant 1,$
\begin{eqnarray}
\lefteqn{
V_1^*(x)\!=\!\frac{KM}{\sqrt{|x|}}e\left(\frac{-x}{2\pi}+\frac{(K-1)^2}{4\pi x}\right)\!
\sum\limits_{l=0}^{L_1}\sum\limits_{0\leqslant l_1\leqslant 2l}
 \sum\limits_{\frac{l_1}{4}\leqslant l_2\leqslant L_2}\!\!\!\!\!\!
c_{l, l_1, l_2}\frac{M^{2l-l_1}(K-1)^{4l_2-l_1}
}{x^{l+3l_2-l_1}}
 }\nonumber\\
&&\hspace{1.5cm}\times 
\left[{\hat{u}_0}^{(2l-l_1)}\left(\frac{(K-1)M}{2\pi x}\right)
+\frac{4\pi^6ix}{45M^6}(t^6\hat{u}_0(t))^{(2l-l_1)}\left(\frac{(K-1)M}{2\pi x}\right)\right]\nonumber\\
&&\hspace{3cm}+O\left(\frac{KM}{\sqrt{|x|}}\left(\frac{T^4}{|x|^3}\right)^{L_2+1}
+K\left(\frac{M}{\sqrt{|x|}}\right)^{2L_1+3}+\frac{K|x|}{M^8}\right)\nonumber
\end{eqnarray}
where $c_{l, l_1, l_2}$ are constants depending only on $l, l_1$ and $l_2.$ 
\end{proposition}
Now we consider the contribution to (A.5) from $V_2(x)$ given by (A.13).
Set
$$
V_2^*(x)=K\int\limits_{-\infty}^\infty\hat{u}_0(t)e\left(\frac{-(K-1)t}{M}+\frac{x}{2\pi}\sin\frac{2\pi t}{M}\right)dt,
$$
then
$$
V_2(x)=\frac{V_2^*(x)-V_2^*(-x)}{2i}.
$$
One can see that $V_2^*(x)$ and $Y_{m, n}^*(x)$ (see (5.8)) have similar integral representation, so they have similar asymptotic behavior (see Propositon (5.1)):
\begin{proposition}
{\rm 1)} For $|x|\geqslant 100 K$ or $|x|\leqslant \frac{1}{100}K,$
$$
V_2^*(x)\ll K^{-A}
$$
where $A>0$ is arbitrarily large and the implied constant depends only on $A$.\\
{\rm 2)} For $\frac{1}{100}K\leqslant |x|\leqslant 100K$, $K^{\frac{3}{8}+\varepsilon}\leqslant M\leqslant K^{\frac{1}{2}}$ and $L_2\geqslant 1,$
\begin{eqnarray}
\lefteqn{
V_2^*(x)=K\sum\limits_{l=0}^{L_2}\sum\limits_{j=0}^l a_{j, l}
\frac{x^l}{M^{5l-2j}}{u_0}^{(5l-2j)}\left(\frac{x-K+1}{M}\right)}\nonumber\\
&&\hspace{4cm}+O\left(\frac{K|x|^{L_2+1}}{M^{3L_2+3}}+\frac{K|x|}{M^7}\right),
\nonumber
\end{eqnarray}
where $a_{j, l}$ are constants depending only on $j$ and $l$.
\end{proposition}
Replacing $T$ by $(K-1)/2$ and $k^*$ by $u_0$ in sections 4 and 5,  one can see that Theorem A.1 follows directly 
from Propositions A.2 and A.3. \\

\noindent {\bf Acknowledgements}\\

\noindent The author would like to thank Dorian Goldfeld and Peter Sarnak for many important comments on this paper.
This work is based on the author's previous work experience with Professors Dorian Goldfeld, Henryk Iwaniec and Peter Sarnak.
It is her pleasure to acknowledge them for the opportunity to share their ideas. She would also like to thank Professor Wenzhi 
Luo for reading the paper and for his comments.\\

\noindent Department of Mathematics,\\
State University of New York at Buffalo, Buffalo, NY, 14260.\\
Email address: xl29@math.buffalo.edu


\begin{thebibliography}{50}
\bibitem[Bu]{}D. Burgess, {\it{On character sums and L-series II,}} Proc. London Math. Soc. 313
(1963), 24-36.
\bibitem[BR]{}J. Bernstein and A. Reznikov, {\it{Periods, subconvexity of $L$-functions and representation theory,}} 
J. Differential Geom. 70 (2005), no. 1, 129--141.
\bibitem[CI]{}J.B. Conrey and H. Iwaniec, {\it{The cubic moment of central values of automorphic $L$-functions, }}
 Ann. of Math. (2) 151 (2000), no. 3, 1175--1216.
\bibitem[De]{} P. Deligne, {\it{Formes modulaires et reprÃÂ©sentations $l-$adiques,}} Lecture Notes in Mathematics 179, 
Springer-Verlag, Berlin/New York (1971) p. 139--172.
\bibitem[DFI1]{} W. Duke, J. Friedlander and H. Iwaniec, {\it{Bounds for automorphic L-functions,}}
Invent. Math. 112 (1993), 1-18.
\bibitem[DFI2]{} W. Duke, J. Friedlander and H. Iwaniec, {\it{Bounds for automorphic L-functions II,}}
Invent. Math. 115 (1994).
\bibitem[DFI3]{}W. Duke, J. Friedlander and H. Iwaniec, {\it{The subconvexity problem for Artin $L$-functions, }}Invent. Math. 149 (2002), no. 3, 489--577.
\bibitem[Gol]{}D. Goldfeld, {\it{Automorphic Forms and L-Functions for the Group
${\rm GL}(  n, \mathbb{R})$,}} Cambridge Studies in Advanced Mathematics, no. 99, 2006.
\bibitem[Go]{}  A. Good, {\it{The square mean of a Dirichlet series associated with cusp forms,}}
Mathematika 29 (1982), 278-295.
\bibitem[GL]{}D. Goldfeld and X. Li,{\it{Voronoi formulas on ${\rm GL}(n)$,}} Int. Math. Res. Not. 2006, Art. ID 86295, 25 pp.
\bibitem[GR]{}I.S. Gradshteyn and I.M. Ryzhik, {\it{Table of integrals, series, and products,}} 
Translated from the Russian. Sixth edition. Academic Press, 2000.
\bibitem[HM]{} G. Harcos and P. Michel, {\it{The subconvexity problem for Rankin-Selberg $L-$functions and Equidistribution of Heegner points. II, }}
 Inventiones Mathematicae 163, 3 (2006), p. 581-655.
\bibitem[He]{}D.R. Heath-Brown, {\it{ The twelfth power moment of the Riemann-function,}} Quart. J. Math. Oxford Ser. (2) 29 (1978), no. 116, 443--462.
\bibitem[Hu]{}M.N. Huxley, {\it{Area, lattice points, and exponential sums,}} London Mathematical Society Monographs. New Series, 13. Oxford Science Publications. 
The Clarendon Press, Oxford University Press, New York, 1996.
\bibitem[Iv]{} A. Ivi\'c, {\it{ On the ternary additive divisor problem and the sixth moment of the zeta-function.}} Sieve methods, exponential sums, and their applications in number theory 
(Cardiff, 1995), 205--243, London Math. Soc. Lecture Note Ser., 237, Cambridge Univ. Press, Cambridge, 1997.
\bibitem[Iw1]{}H. Iwaniec, {\it{Topics in classical automorphic forms. }}Graduate Studies in Mathematics, 17. American Mathematical Society, Providence, RI, 1997. 
\bibitem[Iw2]{} H. Iwaniec, {\it{ Spectral methods of automorphic forms. }}Second edition. Graduate Studies in Mathematics, 53. American Mathematical Society, Providence, RI; Revista MatemÃÂ¡tica Iberoamericana, Madrid, 2002. 
 \bibitem[IK]{}H. Iwaniec and E. Kowalski, {\it{ Analytic number theory,}} 
American Mathematical Society Colloquium Publications, 53. American Mathematical Society, Providence, RI, 2004.
\bibitem[ILS]{}H. Iwaniec, W. Luo, and P. Sarnak, {\it{Low lying zeros of families of $L$-functions. }} Inst. Hautes Ãtudes Sci. Publ. Math.  No. 91  (2000), 55--131 (2001).
\bibitem[JS]{}H. Jacquet and J. Shalika, {\it{ Exterior square $L$-functions,}} Automorphic forms, Shimura varieties, and $L$-functions, Vol. II (Ann Arbor, MI, 1988),  143--226, Perspect. Math., 11, Academic Press, Boston, MA, 1990. 
\bibitem[La]{}E. Lapid,{\it{
On the nonnegativity of Rankin-Selberg $L$-functions at the center of symmetry,}} Int. Math. Res. Not. 2003, no. 2, 65--75.
\bibitem[LYL]{}Y. Lau, J. Liu and Y. Ye, {\it{ A new bound $k\sp {2/3+\epsilon}$ for Rankin-Selberg $L$-functions for Hecke congruence subgroups,}} IMRP Int. Math. Res. Pap. 2006, Art. ID 35090, 78 pp.
\bibitem[Li]{}X. Li, {\it{The central value of the Rankin-Selberg $L-$functions,}} To appear in GAFA.
\bibitem[LRS]{} W. Luo, Z. Rudnick and P. Sarnak, {\it{ On the generalized Ramanujan conjecture for ${\rm GL}(n)$,}} Automorphic forms, automorphic representations, and arithmetic 
(Fort Worth, TX, 1996), 301--310, Proc. Sympos. Pure Math., 66, Part 2, Amer. Math. Soc., Providence, RI, 1999.
\bibitem[KMV]{} E. Kowalski, P. Michel and J.M. Vanderkam, {\it{On the Rankin-Selberg L functions in the level aspect, }}
Duke Mathematical Journal 114(1) (2002). 123-191.
\bibitem[Me]{}T. Meurman, {\it{ On the order of the Maass $L$-function on the critical line,}}  Number theory, Vol. I (Budapest, 1987),  325--354, Colloq. Math. Soc. JÃÂ¡nos Bolyai, 51, North-Holland, Amsterdam, 1990.
\bibitem[Mi]{}P. Michel, {\it{The subconvexity problem for Rankin-Selberg L functions and Equidistribution of Heegner points, }}
Annals of Mathematics 160.1 (2004), 185-236.
\bibitem[MV1]{} P. Michel and A. Venkatesh, {\it{Heegner points and nonvanishing of Rankin-Selberg L-functions,}}Proceedings of Gauss-Dirichlet conference, Clay math proceedings, volume 7.
\bibitem[MV2]{} P. Michel and A. Venkatesh, {\it{ Equidistribution, $L$-functions and ergodic theory: on some problems of Yu. Linnik, }} International Congress of Mathematicians. Vol. II,  421--457, Eur. Math. Soc., ZÃÂ¼rich, 2006. 
\bibitem[Mi]{} S. D. Miller, {\it{
Cancellation in additively twisted sums on ${\rm GL}(n)$,}} Amer. J. Math. 128 (2006), no. 3, 699--729.
\bibitem[MS]{} S.D. Miller and W. Schmid, {\it{ Automorphic distributions, $L$-functions, 
and Voronoi summation for ${\rm GL}(3)$,}} Ann. of Math. (2) 164 (2006), no. 2, 423--488.
\bibitem[Sa]{} P. Sarnak, {\it{Estimates for Rankin-Selberg $L$-functions and quantum unique ergodicity,}} J. Funct. Anal. 184 (2001), no. 2, 419--453.
\bibitem[So]{}D. Soudry, {\it{On Langlands functoriality from classical groups to ${\rm GL}\sb n$. }}Automorphic forms. I.  AstÃÂ©risque  No. 298  (2005), 335--390.
\bibitem[Ti]{} E.C. Titchmarsh, {\it{The theory of the Riemann zeta-function,}} Second edition. 
Edited and with a preface by D. R. Heath-Brown. The Clarendon Press, Oxford University Press, New York, 1986. 
\bibitem[Ve]{}  A. Venkatesh, {\it{Sparse equidistribution problems, period bounds, and subconvexity, }} to appear in Ann. of Math.
\bibitem[Wa]{} G. Watson, {\it{A treatise on the theory of Bessel functions,}} Reprint of the second (1944) edition. Cambridge Mathematical Library. Cambridge University Press, Cambridge, 1995. viii+804 pp.
\bibitem[We]{}H. Weyl, {\it{Zur Absch$\ddot{a}$tzung von $\zeta (1+it)$,}}  Math. Z. 10 (1921), 88-101.
\end{thebibliography}
\end{document}